\documentclass[11pt,a4paper]{article}
\usepackage{amsthm, amssymb, amsmath, a4wide}
\usepackage[all]{xy} 
\usepackage{paralist} 
\usepackage[sans]{dsfont} 
\usepackage{varioref} 
\usepackage{graphicx} 
\usepackage{epsfig}
\usepackage{hyperref}
\usepackage{multicol}

\usepackage{tikz}
\usepackage{pgf}
\usetikzlibrary{calc} 

\usepackage{url} 
\makeatletter

\newcommand{\Rmnum}[1]{\expandafter\@slowromancap\romannumeral #1@}
\makeatother

\let\oldref\ref

\DeclareMathOperator\Div{div} 
 
\DeclareMathOperator\gr{gr}
 
\DeclareMathOperator\ord{ord} 
 
\DeclareMathOperator\pole{pole}
\DeclareMathOperator\proj{Proj}

\DeclareMathOperator\sing{Sing} 
\DeclareMathOperator\spec{Spec}
\DeclareMathOperator\supp{Supp}

\newcommand{\scrA}{\ensuremath{\mathcal{A}}}

\newcommand{\scrC}{\ensuremath{\mathcal{C}}}
\newcommand{\scrD}{\ensuremath{\mathcal{D}}}

\newcommand{\scrF}{\ensuremath{\mathcal{F}}}

\newcommand{\scrI}{\ensuremath{\mathcal{I}}}

\newcommand{\scrL}{\ensuremath{\mathcal{L}}}

\newcommand{\scrP}{\ensuremath{\mathcal{P}}}
\newcommand{\scrQ}{\ensuremath{\mathcal{Q}}}

\newcommand{\cc}{\ensuremath{\mathbb{C}}}
\newcommand{\kk}{\ensuremath{\mathbb{K}}}
\newcommand{\nn}{\ensuremath{\mathbb{N}}}
\newcommand{\pp}{\ensuremath{\mathbb{P}}}
\newcommand{\qq}{\ensuremath{\mathbb{Q}}}
\newcommand{\rr}{\ensuremath{\mathbb{R}}}
\newcommand{\zz}{\ensuremath{\mathbb{Z}}}

\newcommand{\affine}[2]{\ensuremath{\mathbb{A}^{#1}(#2)}}
\newcommand{\projective}[2]{\ensuremath{\pp^{#1}(#2)}}

\newcommand{\ank}{\affine{n}{\kk}}

\newcommand{\pNk}{\projective{N}{\kk}}
\newcommand{\pnk}{\projective{n}{\kk}}

\newcommand{\ooo}{\ensuremath{\mathfrak{o}}}
\newcommand{\ppp}{\ensuremath{\mathfrak{p}}}

\newcommand{\sheaf}{\ensuremath{\mathcal{O}}}

\newcommand{\im}{\ensuremath{\Rightarrow}}

\newcommand{\dsum}{\ensuremath{\bigoplus}}
\newcommand{\into}{\ensuremath{\hookrightarrow}}

\newtheorem{thm}{Theorem}[section]
\newtheorem*{thm*}{Theorem}
\newtheorem{lemma}[thm]{Lemma}
\newtheorem*{lemma*}{Lemma}
\newtheorem{lemma-in-thm}{Lemma}[thm]

\newtheorem{prop}[thm]{Proposition}
\newtheorem*{prop*}{Proposition}
\newtheorem{cor}[thm]{Corollary}
\newtheorem{claim}{Claim}[thm]
\newtheorem*{claim*}{Claim}
\newtheorem{prolemma}[claim]{Lemma}
\newtheorem*{conjecture*}{Conjecture}

\theoremstyle{definition} 
\newtheorem{conjecture}[thm]{Conjecture}
\newtheorem{defn}[thm]{Definition}
\newtheorem*{defn*}{Definition}

\newtheorem*{definotation*}{Definition-Notation}
\newtheorem{example}[thm]{Example}
\newtheorem*{example*}{Example}

\newtheorem*{fact*}{Fact}

\newtheorem*{facts*}{Facts}

\newtheorem{notation}[thm]{Notation}
\newtheorem{note}[thm]{Note}
\newtheorem*{bold-note*}{Note}
\newtheorem{bold-question}[thm]{Question}
\newtheorem{rem}[thm]{Remark}
\newtheorem{reminition}[thm]{Remark-Definition}
\newtheorem*{reminition*}{Remark-Definition}

\newtheorem*{remexample*}{Remark-Example}

\newtheorem*{remtation*}{Remark-Notation}
\newtheorem{remuestion}[thm]{Remark-Question}
\newtheorem*{remuestion*}{Remark-Question}

\newtheorem*{constrinition*}{Construction-Definition}

\theoremstyle{remark}
\newtheorem*{rem*}{Remark}
\newtheorem*{note*}{Note}
\newtheorem*{notation*}{Notation}
\newtheorem*{question*}{Question}

\newtheorem*{questions*}{Questions}

\theoremstyle{plain}

\newcounter{Cases}

\newcounter{UnorderedProofTempCtr}
\newcommand{\tempcommand}{} 
\newenvironment{unorderedproof}[3][]{
	\setcounter{UnorderedProofTempCtr}{\value{#2}}	
	\setcounter{#2}{#3}
	\renewcommand{\tempcommand}{#2}
	\begin{proof}[#1]
} 
{
	\setcounter{\tempcommand}{\value{UnorderedProofTempCtr}}
	\end{proof}
}

\newcommand{\nktorus}{(\kk^*)^n}

\newcommand{\gff}{semidegree} 
\newcommand{\Gf}{\Gff\ }
\newcommand{\Gff}{Semidegree}
\newcommand{\sgff}{subdegree} 
\newcommand{\Sgf}{\Sgff\ }
\newcommand{\Sgff}{Subdegree}

\newcommand{\genfiltrn}{\genfiltrnn\ }
\newcommand{\genfiltrnn}{filtration}

\newcommand{\gengf}{\gengff\ }
\newcommand{\gengff}{\gff}

\newcommand{\gensgf}{\gensgff\ }
\newcommand{\gensgff}{\sgff}

\newcommand{\filtrationchar}{\ensuremath{\mathcal{F}}}
\newcommand{\filtrationring}{\ensuremath{A}}
\newcommand{\profing}{\profingg{\filtrationring}{\filtrationchar}}
\newcommand{\profingg}[2]{\ensuremath{{#1}^{#2}}}
\newcommand{\profinggg}[1]{\profingg{\filtrationring}{#1}}

\newcommand{\gring}{\gringg{\filtrationring}{\filtrationchar}}
\newcommand{\gringg}[2]{\ensuremath{\gr {#1}^{#2}}}

\DeclareMathOperator{\NE}{NE}
\DeclareMathOperator{\nef}{Nef}

\renewcommand{\affine}[2]{\ensuremath{{#2}^{#1}}}

\newcommand{\xdelta}{\bar X^\delta}
\newcommand{\xtildedelta}{\bar X^{\tilde \delta}}
\newcommand{\xf}{\bar X^\scrF}

\begin{document}

\title{Projective completions of affine varieties via degree-like functions}
\author{Pinaki Mondal}
\date{}
\maketitle

\begin{abstract}
We study projective completions of affine algebraic varieties induced by filtrations on their coordinate rings. In particular, we study the effect of the `multiplicative' property of filtrations on the corresponding completions and introduce a class of projective completions (of arbitrary affine varieties) which generalizes the construction of toric varieties from convex rational polytopes. As an application we recover (and generalize to varieties over algebraically closed fields of arbitrary characteristics) a `finiteness' property of divisorial valuations over complex affine varieties proved in \cite{fernex-ein-ishii}. We also find a formula for the pull-back of the `divisor at infinity' and apply it to compute the matrix of intersection numbers of the curves at infinity on a class of compactifications of certain affine surfaces. 
\end{abstract}

\section{Introduction} \label{sec-intro}
\subsection{Summary}
It is a standard fact in algebraic geometry that filtrations (satisfying some natural conditions) on the ring of regular functions of affine varieties correspond to their projective completions (i.e.\ complete varieties containing them as dense open subsets). In this article we study the connection between properties of filtrations and corresponding completions. We primarily focus on (completions corresponding to) filtrations induced by {\em subdegrees} - a notion motivated by the theory of toric varieties (in particular completions of the $n$-torus $(\cc^*)^n$ induced by convex rational polytopes). In short the main results are as follows:\\

\begin{asparaenum}
\item We develop a structure theorem for completions determined by subdegrees and show that in a (precise) sense the correspondence between completions determined by general filtrations and completions determined by subdegrees is analogous to the correspondence between varieties and their normalizations. \\

\item As an application we show that every divisorial valuations over an affine variety is `determined' by its values on a finite number of functions. This result was proved for varieties over $\cc$ in \cite{fernex-ein-ishii} using desingularization and geometry of {\em arc spaces}. Our proof uses essentially elementary (but perhaps somewhat involved) commutative algebra and algebraic geometry. \\

\item We prove a formula for the pull-back (under a dominating morphism) of the `divisor at infinity' on completions determined by subdegrees in terms of the {\em linking numbers} (introduced by Samuel \cite{samuel}) of valuations centered at infinity; in particular this gives a {\em global interpretation} of linking numbers. As applications we compute the matrix of intersection numbers of curves at infinity on a class of completions of certain affine surfaces and find a necessary condition for the finite generation of a subdegree on the coordinate ring of these surfaces. In particular it turns out that the {\em matrix of linking numbers} of the valuations centered at infinity is invertible, which motivates a somewhat amusing conjecture on the invertibility of the matrix of the maximum ratios of a finite collection of $n$-tuples of positive real numbers. \\
\end{asparaenum}

\subsection{More Elaborate Discussion of Main Results} \label{sec-intraborate}

\begin{notation}
Throughout this article $X$ denotes an affine variety over an algebraically closed field $\kk$. The ring of regular functions (resp.\ the field of rational functions) on $X$ will be denoted by $\kk[X]$ (resp.\ $\kk(X)$). Given an ideal $I$ of a ring $A$ (resp.\ a graded ring $B$), the zero set of $I$ in $\spec A$ (resp. $\proj B$) will be denoted by $V(I)$. 
\end{notation}

\paragraph{Basic Construction and Guiding Examples:} We consider the following `rough' correspondence (here we avoid some technical assumptions - see Proposition \ref{prop-spec-completion} for the formal version) - \\

\begin{tabular}{ccc}
\parbox{.25\textwidth}{
Filtrations on $\kk[X]$
} & 
$\longleftrightarrow$ & 
\parbox{.6\textwidth}{
Closed immersions of $X$ into weighted projective spaces. 
} 
\end{tabular}\\

In our case a filtration $\scrF$ on the ring $\kk[X]$ of regular functions on $X$ is a family $\{F_k: k \in \zz\}$ of vector subspaces of the ring $\kk[X]$ of regular functions on $X$ such that $1 \in F_0 \subseteq F_1 \subseteq F_2 \subseteq \cdots \subseteq  \bigcup_{d \geq 0} F_d = \kk[X]$, and $F_dF_e \subseteq F_{d+e}$ for all $d, e \in \zz$. The above correspondence is induced via the natural open immersion $X \into \xf := \proj \left( \sum_{d \geq 0} F_dt^d \right)$, where $t$ is an indeterminate over $\kk[X]$ (Proposition \ref{prop-spec-completion}). The complement of $X$ in $\xf$ is the zero set of the ideal of $t$. In most classical settings filtrations appear in the guise of {\em degree-like functions}, i.e.\ functions $\delta: \kk[X] \to \zz$ which satisfy:

\begin{enumerate}
\item \label{deg1} $\delta(f+g) \leq \max\{\delta(f), \delta(g)\}$ for all $f, g \in \kk[X]$, with $<$ in the preceding equation implying $\delta(f) = \delta(g)$.
\item \label{deg2} $\delta(fg) \leq \delta(f) + \delta(g)$ for all $f, g \in \kk[X]$.
\end{enumerate}

The filtration $\scrF^\delta := \{F^\delta_d\}_{d \in \zz}$ induced by $\delta$ is given by $F^\delta_d := \{f \in \kk[X]: \delta(f) \leq d\}$ and the completion of $X$ corresponding to $\delta$ is $\xdelta := \proj \kk[X]^\delta$, where $\kk[X]^\delta := \sum_{d \geq 0}F^\delta_dt^d$ is the graded ring associated to the filtration constructed as in the preceding paragraph. For example, if $X := \ank$ with coordinates $(x_1, \ldots, x_n)$ and $\delta$ is the degree of polynomials in $\kk[x_1, \ldots, x_n]$, then $\xdelta \cong \pnk$. Similarly, if $\delta$ is the weighted degree on $\kk[x_1, \ldots,x_n]$ which assigns a positive integral weight $d_j$ to $x_j$, $1 \leq j \leq n$, then $\xdelta$ is the weighted projective space $\pp^n(\kk; 1, d_1, \ldots, d_n)$.    \\

A bigger class of completions induced by degree-like functions are projective {\em toric varieties}. Indeed, let $\scrP$ be an $n$-dimensional convex polytope in $\rr^n$ which contains the origin in the interior and whose vertices have rational coordinates. Then $\scrP$ corresponds to an $n$-dimensional projective normal toric variety $X_\scrP$. We claim that $X_\scrP$ is the completion of the $n$-torus $X := (\kk^*)^n$ induced by a degree-like function. Indeed, define $\delta_\scrP: \kk[x_1, x_1^{-1}, \ldots, x_n, x_n^{-1}] \setminus\{0\} \to \qq$ such that 
$$\delta_\scrP(\sum a_\alpha x^\alpha) := \inf\{r \in \rr: r \geq 0,\ \alpha \in r\scrP\ \text{for all}\ \alpha \in \zz^n\ \text{such that}\ a_\alpha \neq 0\}$$ 
(see Figure \ref{sub-picture} for an example). Pick any positive integer $e$ such that $e\delta_\scrP$ is integer-valued. Then $e\delta_\scrP$ is a degree-like function on $\kk[X]$ and it is not hard to see that $X_\scrP \cong \bar X^{e\delta_\scrP}$. \\
\begin{figure}[h]
\begin{center}
\begin{tikzpicture}[scale=0.6]
	\begin{scope}[shift={(0,0)}]
		\draw [gray,  line width=0pt] (-1.5,-1.5) grid (2.5,2.5);
		\draw [<->] (0,2.5) node (yaxis) [above] {$y$}
       	 |- (2.5,0) node (xaxis) [right] {$x$};
       	 
       	\draw[green,thick ] (-1,-1) -- (2,-1) -- (-1,2) -- cycle;
       	\fill[red] (0,0) circle (4pt);
       	
       	\draw (-0.5,-0.5) node {\textcolor{green}{$\scrP$}};
		\draw (-2.5,-1) node [below right, text width= 4.5cm] {
			\scriptsize{
 			\begin{align*}
 			\delta_{\scrP}: &~ x \mapsto 1, x^{-1} \mapsto 1, y \mapsto 1 \\
 							&~ y^{-1} \mapsto 1, x^{-1}y^{-1} \mapsto 1 \\
 							&~ x^2y^{-1} \mapsto 1, x^{-1}y^2 \mapsto 1
 			\end{align*}
 			}
		};    	

	\end{scope}
	
	\begin{scope}[shift={(7,0)}]
		\draw [gray,  line width=0pt] (-1.5,-1.5) grid (2.5,2.5);
		\draw [<->] (0,2.5) node (yaxis) [above] {$y$}
       	 |- (2.5,0) node (xaxis) [right] {$x$};
       	 
       	 \fill[red] (0,0) circle (4pt);
       	 
       	\draw[green,thick ] (-1,-1) -- (2,-1);

		\draw [blue, thick, ->] (0,-1) -- (0,-1.5);
		
		\draw (-2.5,-1) node [below right, text width= 4cm] {
			\scriptsize{
 			\begin{align*}
 				\delta_1(x^\alpha y^\beta) 	&= (\alpha, \beta) \cdot (0,-1) \\
 											&= -\beta
 			\end{align*}
 			}
		};   
	\end{scope}
	
	\begin{scope}[shift={(13,0)}]
		\draw [gray,  line width=0pt] (-1.5,-1.5) grid (2.5,2.5);
		\draw [<->] (0,2.5) node (yaxis) [above] {$y$}
       	 |- (2.5,0) node (xaxis) [right] {$x$};
       	 
       	 \fill[red] (0,0) circle (4pt);
       	 
       	\draw[green,thick ] (2,-1) coordinate (p1) -- (-1,2) coordinate (p2);
       	\coordinate (m) at ($(p1)!0.5!(p2)$);
		\draw [blue, thick, ->] (m) -- ($(m)!0.5cm!-90:(p2)$);
		
		\draw (-2.5,-1) node [below right, text width= 4cm] {
			\scriptsize{
 			\begin{align*}
 				\delta_2(x^\alpha y^\beta) 	&= (\alpha, \beta) \cdot (1,1) \\
 											&= \alpha + \beta
 			\end{align*}
 			}
		};

	\end{scope}
	
	\begin{scope}[shift={(19,0)}]
		\draw [gray,  line width=0pt] (-1.5,-1.5) grid (2.5,2.5);
		\draw [<->] (0,2.5) node (yaxis) [above] {$y$}
       	 |- (2.5,0) node (xaxis) [right] {$x$};
       	 
       	 \fill[red] (0,0) circle (4pt);
       	 
       	\draw[green,thick ] (-1,-1) -- (-1,2);

		\draw [blue, thick, ->] (-1,0) -- (-1.5,0);
		
		\draw (-2.5,-1) node [below right, text width= 4cm] {
			\scriptsize{
 			\begin{align*}
 				\delta_3(x^\alpha y^\beta) 	&= (\alpha, \beta) \cdot (-1,0) \\
 											&= -\alpha
 			\end{align*}
 			}
		};   
	\end{scope}
\end{tikzpicture}
\caption{$\delta_\scrP = \max\{\delta_1, \delta_2, \delta_3\}$}  \label{sub-picture}
\end{center}
\end{figure}

Note that our first two examples of degree-like functions, namely degree and weighted degrees on polynomial rings, satisfy the multiplicative property (i.e.\ property \ref{deg2}) of degree-like functions with exact equality. On the other hand, it is straightforward to see that even though the degree-like functions $\delta_\scrP$ of the preceding paragraph in general do {\em not} satisfy property \ref{deg2}, they are the maxima of finitely many degree-like functions which {\em do} satisfy it, namely the weighted degrees determined by the facets (i.e.\ $(n-1)$-dimensional faces) of $\scrP$ (see Figure \ref{sub-picture} and Example \ref{toric-example}). Motivated by these examples we define two special classes of degree-like functions: {\em semidegrees}, which satisfy property \ref{deg2} of degree-like functions with exact equality, and {\em subdegrees}, which are the maxima of finitely many semidegrees. The primary motive of this article is to (start to) explore to what extent properties of toric varieties remain true for completions of affine varieties determined by general subdegrees. In a sequel to this paper we plan to study compactifications of the affine space determined by special classes of subdegrees, e.g.\ subdegrees determined by weighted degrees in {\em different} systems of coordinates. We would like to mention that when dealing with compactifications, sometimes using {\em only} the formalism of degree-like functions may simplify problems, as we found out in the proof of the `algebraicity criterion' of \cite{contractibility}.

\begin{rem} \label{valuative-remark}
Note that the negative of an (integer-valued) semidegree is simply a {\em discrete valuation}. More generally, negative of degree-like functions have been studied, from the perspective of {\em local} geometry or algebra. 
\begin{tabbing}
Degree-like function	\= blaank 	\= Classical terminologies for its negative blah blah blah \= \kill
{\bf Our notion}		\> 				\> \bf{Classical terminologies for its negative}\\
degree-like function	\> 				\> order function \cite{szpiro}, \cite{lejier}, pseudo-valuation \cite{huckaba} \\
subdegree				\>				\> homogeneous order function \cite{szpiro}, subvaluation \cite{huckaba}\\
semidegree				\>				\> (discrete) valuation
\end{tabbing}
Our `structure theorem' for subdegrees (Theorem \ref{gengfsgf-characterization}) is a strengthening of a special case of \cite[Theorem 1]{szpiro} (see Remark \ref{szpiro-sgf}). An essential ingredient of our `main existence theorem' (Theorem \ref{thm-noeth-integral-subdegree}) is Rees' theory of subvaluations associated with ideals \cite{rees-vals-ideals-i}, \cite{rees-vals-ideals-ii}.
\end{rem}

\paragraph{Structure Theorem for Subdegrees and its Corollaries:} Let $\delta$ be a subdegree on $\kk[X]$ with $\delta = \max\{\delta_1, \ldots, \delta_k\}$. Getting rid of the `extra' $\delta_j$'s if necessary, we may assume the presentation is {\em minimal}, i.e.\ for every $j$, $1 \leq j \leq k$, there exists $f \in \kk[X]$ such that $\delta_j(f) > \delta_i(f)$ for all $i \neq j$, $1 \leq i \leq k$. Our first main result is the `structure theorem' for subdegrees:

\begin{thm}[Informal version of Theorem \ref{gengfsgf-characterization} and Proposition \ref{pole-and-degree}] \label{weak-structure}
Every subdegree $\delta$ has a unique minimal presentation. A degree-like function $\delta$ is a subdegree iff the ideal $I$ of $\kk[X]^\delta$ generated by $t$ is a decomposable radical ideal. If $\delta$ is a subdegree then the non-zero semidegrees in the minimal presentation of $\delta$ are multiples of the orders of pole (of rational functions on $X$) along the irreducible components of $V(I) = \xdelta \setminus X$. 
\end{thm}

An immediate corollary of Theorem \ref{weak-structure} gives a characterization of subdegrees in the case that $\kk[X]^\delta$ is Noetherian, namely: $\delta$ is a subdegree iff $\delta(f^k) = k\delta(f)$ for all $f \in \kk[X]$ and $k \geq 0$. Another corollary of the arguments in the proof of Theorem \ref{weak-structure} is a `finiteness property' of {\em divisorial valuations}. Recall that divisorial valuations on $\kk(X)$ are valuations of the form $\nu := q \ord_E : \kk(X)^* \to \zz$, where $E$ is a prime divisor on a normal variety $Y$ equipped with a proper birational morphism $f : Y \to X$, $q = q(v)$ is a positive integer number called the {\em multiplicity} of $\nu$, and for every $h \in \kk(X)^*$ which is regular at the generic point of $f(E)$, $\ord_E(h)$ is the order of vanishing of $h \circ f$ at the generic point of $E$. Given a divisorial valuation $\nu$ over $X$, we may consider that $\nu$ is {\em centered at infinity} with respect to some affine variety $Z$ which is birational to $X$ and therefore we may assume by means of Theorem \ref{weak-structure} that $-\nu$ is a semidegree associated to some subdegree on $\kk[Z]$. Applying arguments of the proof of Theorem \ref{weak-structure} to this context yields the following result which was proved in \cite[Theorem 0.2]{fernex-ein-ishii} for $\kk = \cc$.

\begin{thm*}[Theorem \ref{finiteness-prop}]
Let $X$ be an irreducible affine variety and $\nu$ be a divisorial valuation over $X$. Then there exist elements $f_1, \ldots, f_r \in \kk[X]\setminus\{0\}$ such that for every $f \in \kk[X] \setminus\{0\}$,
\begin{align}
\nu(f) &= \min\{\nu'(f) : \nu'\ \text{is a divisorial valuation over $X$,}\ \nu'(f_i) = \nu(f_i)\  \text{for}\ 1 \leq i \leq r \} \tag{$*$} \label{finiteness}\\
	&= \min\{\nu'(f) : \nu'\ \text{is a valuation on $\kk(X)$ such that the value group of $\nu$ contains} \notag \\
	& \mbox{\phantom{=$\min\{\nu'(f):$}}\ \ \text{the integers and}\ \nu'(f_i) = \nu(f_i)\ \text{for}\ 1 \leq i \leq r \} \tag{$*'$} \label{finiteness'}
\end{align}
\end{thm*}

\paragraph{Normalization at infinity of completions:} In order to study a given completion $\bar X$ of an affine variety $X$, it is natural to try to find a `simpler' model, i.e.\ a `simpler' completion $\bar X'$ of $X$ which {\em dominates $\bar X$}, i.e.\ which comes with a morphism $\pi: \bar X' \to \bar X$ that restricts to the identity on $X$. If $X$ is normal or non-singular, a natural choice for the `simpler' model is the normalization or a desingularization of $\bar X$ (provided the latter exists of course). For a general $X$, we introduce the {\em normalization at infinity of $\bar X$ with respect to $X$}, which is the unique {\em minimal} completion $\bar X'$ dominating $\bar X$ such that $\bar X'$ is {\em normal at infinity with respect to $X$}, i.e.\ for all open subset $U$ of $\bar X'$, the ring of regular functions on $U$ is integrally closed in the ring of regular functions on $U \cap X$. In particular, if $\bar X$ is normal, then $\bar X' \cong \bar X$. Using the theory of {\em Rees valuations} (\cite{rees-vals-ideals-i}, \cite{rees-vals-ideals-ii}) we show that

\begin{thm*}[A reformulation of Theorem \ref{thm-noeth-integral-subdegree}] 
Let $\bar X$ be a projective completion of $X$ such that $\bar X \setminus X$ is the support of an effective ample divisor. Then there is a subdegree $\delta'$ on $\kk[X]$ such that $\bar X^{\delta'}$ is the normalization at infinity of $\bar X$ with respect to $X$. If $\delta$ is any degree-like function such that $\xdelta \cong \bar X$, then one such $\delta'$ is given by:
$$\delta'(f) := e \lim_{n \to \infty}\frac{\delta(f^n)}{n}\ \text{for all}\ f \in \kk[X],$$
where $e$ is a sufficiently divisible integer. 
\end{thm*}

\begin{note}
A completion $\bar X$ of $X$ is induced by a degree-like function iff $\bar X \setminus X$ is the support of an effective ample divisor (Proposition \ref{prop-spec-completion}).
\end{note}

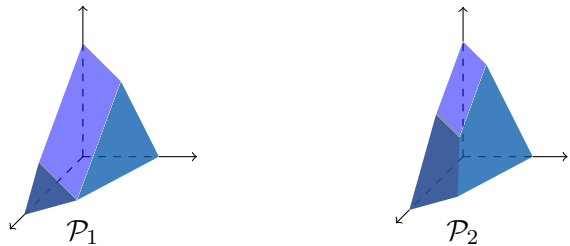
\begin{figure}[h]
\begin{center}
	
\begin{tikzpicture}[scale=1/2]
	\begin{scope}[shift={(0,0)}, scale=36]
    	\def\pointA{\pgfpointxyz{0}{0}{1/9}}
    	\def\pointB{\pgfpointxyz{1/18}{0}{0}}
    	\def\pointC{\pgfpointxyz{0}{1/12}{0}}
    	\def\pointD{\pgfpointxyz{1/36}{0}{1/12}}
    	\def\pointE{\pgfpointxyz{1/36}{1/18}{0}}
    	\def\pointF{\pgfpointxyz{0}{1/36}{1/12}}
    	
    	\def\pointZ{\pgfpointxyz{0}{0}{1/36}}
    	\def\pointX{\pgfpointxyz{1/36}{0}{0}}
    	\def\pointY{\pgfpointxyz{0}{1/36}{0}}
    	
		\pgfsetdash{{0.1cm}{0.1cm}}{0pt}		
		\pgfpathmoveto{\pgfpointorigin}
		\pgfpathlineto{\pointA}
		\pgfpathmoveto{\pgfpointorigin}
		\pgfpathlineto{\pointB}
		\pgfpathmoveto{\pgfpointorigin}
		\pgfpathlineto{\pointC}
		\pgfusepath{stroke}
		
		\pgfsetdash{}{0pt}
		\pgfsetarrowsend{to}
		
		\pgfpathmoveto{\pointA}
		\pgfpathlineto{\pgfpointadd{\pointA}{\pointZ}}
		\pgfusepath{stroke}
		\pgfpathmoveto{\pointB}
		\pgfpathlineto{\pgfpointadd{\pointB}{\pointX}}
		\pgfusepath{stroke}
		\pgfpathmoveto{\pointC}
		\pgfpathlineto{\pgfpointadd{\pointC}{\pointY}}
		\pgfusepath{stroke}		
		
		\color{red}
		\pgfsetfillopacity{0.5}

		\pgfpathmoveto{\pointF}
		\pgfpathlineto{\pointA}
		\pgfpathlineto{\pointD}
		\pgfpathlineto{\pointF}
		\pgfusepath{fill}
		
		\color{green}
		\pgfsetfillopacity{0.5}		
		\pgfpathmoveto{\pointD}
		\pgfpathlineto{\pointB}
		\pgfpathlineto{\pointE}
		\pgfpathlineto{\pointD}
		\pgfusepath{fill}
		
		\color{blue}
		\pgfsetfillopacity{0.5}
		\pgfpathmoveto{\pointD}
		\pgfpathlineto{\pointF}
		\pgfpathlineto{\pointC}
		\pgfpathlineto{\pointE}		
		\pgfpathlineto{\pointD}
		\pgfusepath{fill}

    \end{scope}
    
    \begin{scope}[shift={(10,0)}, scale=55]
    	\def\pointA{\pgfpointxyz{0}{0}{1/15}}
    	\def\pointB{\pgfpointxyz{1/30}{0}{0}}
    	\def\pointC{\pgfpointxyz{0}{1/18}{0}}
    	\def\pointD{\pgfpointxyz{1/60}{0}{1/20}}
    	\def\pointE{\pgfpointxyz{1/90}{4/90}{0}}
    	\def\pointF{\pgfpointxyz{0}{1/30}{1/30}}
    	\def\pointG{\pgfpointxyz{1/90}{1/45}{1/30}}
    	
    	\def\pointZ{\pgfpointxyz{0}{0}{1/60}}
    	\def\pointX{\pgfpointxyz{1/60}{0}{0}}
    	\def\pointY{\pgfpointxyz{0}{1/60}{0}}
    	
    	\pgfsetdash{{0.1cm}{0.1cm}}{0pt}		
		\pgfpathmoveto{\pgfpointorigin}
		\pgfpathlineto{\pointA}
		\pgfpathmoveto{\pgfpointorigin}
		\pgfpathlineto{\pointB}
		\pgfpathmoveto{\pgfpointorigin}
		\pgfpathlineto{\pointC}
		\pgfusepath{stroke}
		
		\pgfsetdash{}{0pt}
		\pgfsetarrowsend{to}
		
		\pgfpathmoveto{\pointA}
		\pgfpathlineto{\pgfpointadd{\pointA}{\pointZ}}
		\pgfusepath{stroke}
		\pgfpathmoveto{\pointB}
		\pgfpathlineto{\pgfpointadd{\pointB}{\pointX}}
		\pgfusepath{stroke}
		\pgfpathmoveto{\pointC}
		\pgfpathlineto{\pgfpointadd{\pointC}{\pointY}}
		\pgfusepath{stroke}		
		
		\color{red}
		\pgfsetfillopacity{0.5}

		\pgfpathmoveto{\pointF}
		\pgfpathlineto{\pointA}
		\pgfpathlineto{\pointD}
		\pgfpathlineto{\pointG}
		\pgfpathlineto{\pointF}
		\pgfusepath{fill}
		
		\color{green}
		\pgfsetfillopacity{0.5}		
		\pgfpathmoveto{\pointD}
		\pgfpathlineto{\pointB}
		\pgfpathlineto{\pointE}
		\pgfpathlineto{\pointG}
		\pgfpathlineto{\pointD}
		\pgfusepath{fill}
		
		\color{blue}
		\pgfsetfillopacity{0.5}
		\pgfpathmoveto{\pointG}
		\pgfpathlineto{\pointF}
		\pgfpathlineto{\pointC}
		\pgfpathlineto{\pointE}		
		\pgfpathlineto{\pointG}
		\pgfusepath{fill}
		
    \end{scope}
    
   	\node at (0,-2){$\scrP_1$};
	\node at (10,-2){$\scrP_2$};
    
\end{tikzpicture}
\caption{Polytopes corresponding to toric varieties giving negative answer to Question \ref{isomorphic-question}}  \label{semi-picture}
\end{center}
\end{figure}

\paragraph{Divisor at Infinity and its Multiplicities:} Let $\delta$ be a subdegree and $\delta = \max\{\delta_1, \ldots, \delta_k\}$ be the minimal presentation of $\delta$ as the maximum of semidegrees. It follows from Theorem \ref{weak-structure} that $\delta_j = -d_j\ord_{E_j}$, where $E_1, \ldots, E_k$ are irreducible components of $\xdelta \setminus X$, and for each $j$, $\ord_{E_j}$ is the order of vanishing along $E_j$ and $d_j$ is some positive integer. This motivates a natural question about the roles of the numbers $d_j$. More precisely, it is natural to ask
\begin{bold-question} \label{isomorphic-question}
If $\bar X_1$ and $\bar X_2$ are two completions of $X$ such that the set of divisorial valuations associated to the irreducible complements of $\bar X_j \setminus X$ are identical, then is it true that $\bar X_1 \cong \bar X_2$?'
\end{bold-question}
If $\dim X \leq 2$, then it is easy to see that the answer to Question \ref{isomorphic-question} is positive. However, straightforward examples from toric geometry shows that it is not true in higher dimensions. Indeed, for $n \geq 3$, it is easy to find two convex rational polytopes (`rational' means that the vertices have rational coordinates) $\scrP_1, \scrP_2$ in $\rr^n$ such that the following properties are satisfied (see Figure \ref{semi-picture}):
\begin{compactenum}
\item there is a one-to-one correspondence between the facets of $\scrP_i$'s such that for each facet $\scrQ_1$ of $\scrP_1$, the corresponding facet $\scrQ_2$ of $\scrP_2$ has the {\em same} outer normal (i.e.\ $\scrQ_1$ and $\scrQ_2$ are parallel and on the `same side' of the origin), but
\item $\scrP_1$ and $\scrP_2$ are {\em not} combinatorially isomorphic (i.e.\ the {\em face lattice} of $\scrP_i$'s are not isomorphic).  
\end{compactenum}

Then it follows from basic theory of toric varieties (e.g.\ as in \cite{fultoric}) that the toric varieties $X_{\scrP_1}$ and $X_{\scrP_2}$ are not isomorphic, even though the irreducible components of the complement of the $n$-torus in $X_{\scrP_j}$'s determine the same set of divisorial valuations.\\

The preceding example shows that $q_j$'s play a non-trivial role in the structure of $\xdelta$. One of the places they manifest is in the divisor on $\xdelta$ determined by $t \in \kk[X]^\delta = \sum_{d \geq 0}F^\delta_d t^d$. Indeed, the $\qq$-Cartier divisor of $t$ on $\xdelta$ is precisely $\sum_{j=1}^k \frac{1}{d_j}E_j$ (Lemma \ref{divisor-at-infinity}). We call this divisor the {\em divisor at infinity} and denote it by $D^\delta_\infty$. The toric analogue for the divisor at infinity is the {\em divisor of a polytope} (\cite[Section 7.1]{littlehalcox}). We prove a formula for its pull-back under a dominant morphism:

\begin{prop*}[Proposition \ref{pulling-back-infinity}] \label{pulling-back-infinity-intro}
Let $\phi: Z \to \xdelta$ be a dominant morphism. Then
\begin{align}
\phi^*(D^\delta_\infty) &= \sum_W l^\phi_\infty(\delta, \pole_W) [W], \label{pull-back-formula-intro}
\end{align}
where the sum is over codimension one irreducible subvarieties $W$ of $Z$, and for each such $W$, the function $\pole_W$ is the negative of the {\em order $\ord_W$ of vanishing} along $W$ and  
$$l^\phi_\infty(\delta, \pole_W)  := \max \left \{\frac{\pole_W(\phi^*(f))}{\delta(f)}: f \in \kk[X],\ \delta(f) > 0\right\}.$$
\end{prop*}   

\begin{rem}
Identity \eqref{pull-back-formula-intro} in particular implies that $l^\phi_\infty(\delta, \pole_W)$ exists. Note also that the sum in identity \eqref{pull-back-formula-intro} is finite, since $l^\phi_\infty(\delta, \pole_W)= 0$ for all codimension one irreducible subvariety $W$ of $Z$ such that $W \cap \phi^{-1}(X) \neq \emptyset$.
\end{rem}

\begin{reminition} \label{linking-reminition}
Given a Krull local ring $\ooo$, a valuation $\omega$ on $\ooo$ and a height one prime ideal $\ppp$ of $\ooo$, the {\em linking number} \cite{huckaba} of $\omega$ and the valuation $\nu_\ppp$ on $\ooo$ corresponding to $\ppp$ is 
$$l(\omega, \nu_\ppp) := \inf_{f \in \ppp,\ f \neq 0}\frac{\omega(f)}{\nu_\ppp(f)}.$$
The linking number was introduced in \cite{samuel} in connection with defining the pull-back of Weil divisors under a birational regular mapping. More precisely, if $\phi: Z \to Y$ is a birational map and $\omega$ (resp.\ $\nu_\ppp$) is the order of vanishing along a codimension one subvariety $W$ of $Z$ (resp.\ $V$ of $Y$), then $l(\omega, \nu_\ppp)$ is a `candidate' for the coefficient of $[W]$ in $\phi^*([V])$. Now assume $Y := \xdelta$ for a {\em semidegree} $\delta$ (i.e.\ $Y\setminus X$ is irreducible) and $V := Y \setminus X$. Then $\delta = -\nu_\ppp$ and $\pole_W = - \omega$. \cite[Theorem 2]{samuel} and identity \eqref{pull-back-formula-intro} then imply that $$l^\phi_\infty(\delta,\pole_W) = l(-\pole_W,-\delta).$$
It is perhaps instructive that $l(-\pole_W,-\delta)$ and $l^\phi_\infty(\delta,\pole_W)$ measure the {\em same} quantity, but the $in\!f$ of the local computation is replaced by a $sup$ in the global computation. We call $l^\phi_{\infty}(\delta, \pole_W)$ the {\em linking number at infinity} ({\em relative to $\phi$}) of $\delta$ and $\pole_W$. In the case that $\phi$ is birational (so that $\phi^*$ is identity on the function field), we simply write $l_{\infty}(\delta, \pole_W)$ for $l^\phi_{\infty}(\delta, \pole_W)$.
\end{reminition}

\paragraph{Intersection numbers of curves at infinity:} Let $X$ be a normal affine surface with no non-constant invertible regular functions (e.g.\ $\affine{2}{\kk}$) and $\bar X_1, \ldots, \bar X_k$ be normal algebraic completions of $X$ such that for each $j$, the complement $C_j$ of $X$ in $\bar X_j$ is an irreducible curve. Let $\bar X$ be the normalization of the diagonal embedding of $X$ into $\bar X_1 \times \cdots \times \bar X_k$. We apply Proposition \ref{pulling-back-infinity} to compute the matrix of intersection numbers of the curves at infinity on $\bar X$ in terms of the linking numbers at infinity of associated semidegrees (Lemma \ref{linking-intersections}) - in particular, it turns out that in this case the matrix of linking numbers at infinity of the semidegrees is {\em invertible}. Taking the semidegrees to be weighted degrees yields the following statement (which is not hard to verify directly, see e.g.\ \cite{mathover-det}):

\begin{cor}[cf.\ Corollary \ref{toric-intersections}] \label{deterally}
Let $n = 2$. Let $v_1, \ldots, v_k$ be pairwise non-proportional elements in $\qq^n$ with positive coordinates and $L$ be the $k \times k$ matrix with entries 
$$l_{ij} := \max\{\frac{v_{jm}}{v_{im}}: 1 \leq m \leq n\}.$$
Then $L$ is invertible.
\end{cor}
Motivated by Corollary \ref{deterally} and a lot of computation (done mostly for $n=3$), we propose the following
\begin{conjecture} \label{deterecture}
The conclusion of Corollary \ref{deterally} remains true for all $n \geq 2$. 
\end{conjecture}
Since $l_{ij}$'s are precisely the linking numbers at infinity of weighted degrees determined by $v_i$ and $v_j$, it seems plausible that Conjecture \ref{deterecture} can be proved/disproved using the geometry of toric varieties associated to polytopes with facets determined by $v_j$'s. 

\paragraph{Finite generation of subdegrees:} For a subdegree $\delta$ on $\kk[X]$ to correspond to a completion of $X$, it is necessary that $\delta$ is {\em finitely generated}, i.e.\ the graded ring $\kk[X]^\delta$ is finitely generated as an algebra over $\kk$. It is therefore natural to try to find conditions which guarantee that $\delta$ is finitely generated. Note that if $\delta_1, \ldots, \delta_k$ are semidegrees such that $\delta = \max\{\delta_1, \ldots, \delta_k\}$, then $\kk[X]^\delta = \bigcap_{j=1}^k \kk[X]^{\delta_j}$. It is therefore natural to ask:
\begin{bold-question} \label{finite-question}
If each $\delta_j$ is finitely generated, then is it true that $\delta$ is also finitely generated?
\end{bold-question}
If $X = (\kk^*)^n$ and each $\delta_j$ is a weighted degree in a fixed set of coordinates on $X$, then the answer to Question \ref{finite-question} is positive and is an implication of {\em Gordan's Lemma} (which states that the integral points in a closed rational cone of $\rr^n$ form a finitely generated semigroup). But as the following example from \cite{subalsection} shows, the answer is in general false, even for $n = 2$.  

\begin{example} \label{non-finite-example}
Consider elements $y_1 := y - x^5 - x^{-2}$ and $y_2 := y + x^5 - x^{-2}$ of $\cc(x,y)$. Let $\delta_1$ (resp.\ $\delta_2$) be the weighted degree on $\cc(x,y)$ with respect to coordinates $(x,y_1)$ (resp.\ $(x,y_2)$) corresponding to weights $1$ for $x$ and $-3$ for $y_1$ (resp.\ $y_2$). Then it turns out that $\cc[x,y]^{\delta_j}$ is finitely generated for each $j$ and $\delta := \max\{\delta_1, \delta_2\}$ is strictly positive on non-constant elements of $\cc[x,y]$, but $\delta$ is not finitely generated - see \cite[Section 2]{subalsection} for details. In fact it is shown in \cite{subalsection} that $\delta$ corresponds to a {\em non-algebraic} normal analytic compactification of $\cc^2$. 
\end{example}

For $n = 2$, we give a positive answer to Question \ref{finite-question} under the additional hypothesis that each $\delta_j$ is strictly positive on every non-constant regular function on $X$ (Corollary \ref{finite-corollary}). Indeed, this additional hypothesis is equivalent to saying that each ${\bar X}^{\delta_j}$ is a projective completion of $X$. Let $\bar X$ be the normalization of the diagonal embedding of $X$ into $\bar X^{\delta_1} \times \cdots \times \bar X^{\delta_k}$. Then Question \ref{finite-question} reduces to a question about ampleness of some divisors on $\bar X$, which we answer using Proposition \ref{pulling-back-infinity}. 

\begin{rem}
Special cases of the surfaces $\bar X$ of the preceding paragraph (for $X = \cc^2$) have been studied in \cite{campillo-piltant-lopez-cones-surfaces}. They constitute a natural class of surfaces for which the cone of curves is `nice', e.g.\ it is finitely generated and simplicial (this follows e.g.\ from Lemma \ref{linking-intersections}).
\end{rem}

\begin{remuestion}
Counter-examples to Question \ref{finite-question} that we know of (e.g.\ Example \ref{non-finite-example}) require the characteristic of the field to be zero. Therefore we ask: is the answer to Question \ref{finite-question} positive if $\kk$ has positive characteristic? 
\end{remuestion}

\subsection{Organization}
In the next section we introduce and give some examples of the correspondence between filtrations on the coordinate ring $\kk[X]$ of an affine variety $X$ and projective completions of $X$. In particular, we show that `finitely generated' filtrations on $\kk[X]$ correspond to immersions of $X$ into {\em generalized weighted projective spaces} (Proposition \ref{prop-spec-completion}). In Section \ref{sec-semiquasintro} we reformulate the notion of filtrations in terms of {\em degree-like functions} and introduce two special classes of degree-like functions which are protagonists of this article, namely {\em semidegree} and {\em subdegree}. \\

The remaining sections contain our main results. In Section \ref{sec-sub-structure} we establish the structure theorem of subdegrees (Theorem \ref{gengfsgf-characterization}) and some of its immediate corollaries. Section \ref{sec-sub-normality} is devoted to the study of `normality' properties of completions determined by subdegrees. In this section we establish the relation between semidegrees in the minimal presentation of a subdegree and orders of pole along components of the hypersurface at infinity on the completion induced by the subdegree (Proposition \ref{pole-and-degree}). We also introduce the notion of {\em relative normality at infinity}, and show that the completion $\xdelta$ of an affine variety $X$ determined by a a subdegree $\delta$ on $\kk[X]$ is relatively normal at infinity with respect to $X$ (Proposition \ref{relatively-normal-subdegree}). We then establish the existence of a subdegree $\tilde \delta$ which {\em normalizes} a given degree-like function $\delta$ in the sense that $\bar X^{\tilde \delta}$ is the normalization at infinity (with respect to $X$) of $\xdelta$ (Theorem \ref{thm-noeth-integral-subdegree}). We end Section \ref{sec-sub-normality} with a proof of the `finiteness property' of divisorial valuations discovered (in the case of complex varieties) in \cite{fernex-ein-ishii} (Theorem \ref{finiteness-prop}). In Section \ref{sec-sub-divisor} we prove the `pull-back formula' of the divisor at infinity (Proposition \ref{pulling-back-infinity}) and apply it to give a geometric proof of Corollary \ref{deterally} and find a sufficient condition for the answer to Question \ref{finite-question} to be positive.\\

This article is mostly based on my PhD thesis. I express my gratitude to my thesis advisor Professor Pierre Milman for posing the questions, helpful suggestions, and in general guiding me throughout this work. I would also like to thank Professor Askold Khovanskii for helpful suggestions (e.g.\ he pointed out the possibility of a connection between semidegrees and `orders of poles at infinity' \cite{khovanication-order}) and Professor Bernard Teissier for bringing the article \cite{szpiro} to my attention.

\section{Compactifications from Filtrations} \label{sec-filtrintro}
\newcommand{\af}{A^\scrF}

In this section we set up the correspondence between filtrations and completions. The material is basic and well-known (except for assertion \ref{compactn-to-filtrn} of Proposition \ref{prop-spec-completion}). 

\begin{defn} A \genfiltrn $\scrF$ on an $\kk$-algebra $A$ is a family $\{F_i: i \in \zz\}$ of $\kk$-vector subspaces of $A$ such that
\begin{compactenum}
\item $F_i \subseteq F_{i+1}$ for all $i \in \zz$ \label{genfiltrn-prop-1}
\item $1 \in F_0\setminus  F_{-1}$ \label{genfiltrn-prop-2}
\item $A = \bigcup_{i \in \zz}F_i$ and \label{genfiltrn-prop-3}
\item $F_iF_j \subseteq F_{i+j}$ for all $i,j \in \zz$. \label{genfiltrn-prop-4}
\end{compactenum}
\end{defn}

\begin{rem}
The condition $1 \not\in F_{-1}$ is introduced only to exclude the trivial filtration, i.e.\ the filtration $\{F_i := A\}_{i \in \zz}$.
\end{rem}

To each filtration $\scrF$ on a $\kk$-algebra $A$ we associate two graded rings:
\begin{gather*}
\af := \dsum_{d \geq 0} F_d,\ \text{and}\ \gring := \dsum_{d \geq 0} (F_d/F_{d-1}).
\end{gather*}
For each $f \in F_d$, we denote the `copy' of $f$ in the $d$-th graded component of $\af$ by $(f)_d$. The ring $\af$ has a natural graded $\kk$-algebra structure with the multiplication defined by: 
$$\left(\sum_d(f_d)_d \right)\left(\sum_e(g_e)_e\right) := \sum_k\sum_{d+e = k}(f_dg_e)_k.$$
We say that $\scrF$ is a {\em finitely generated} \genfiltrn if $\profing$ is a finitely generated $\kk$-algebra.

\begin{rem} \label{t-remark}
Let $t$ be an indeterminate over $A$. Then there is an isomorphism
\begin{align*}
\af \cong \sum_{d \geq 0} F_dt^d \subseteq A[t,t^{-1}]
\end{align*}
which maps $(f)_d \mapsto ft^d$ for all $f \in A$ (in particular, $(1)_1 \mapsto t$). It follows immediately from this isomorphism that $A$ is a domain iff $\profing$ is also a domain.
\end{rem}

\begin{rem}[Motivation for the grading on $\af$]
Note that even though it is possible for a filtration to have non-trivial components with {\em negative} degree, when forming the graded rings we take the direct sum only over components of non-negative degree. The motivation is geometric: given a variety $Y$ and a Cartier divisor $D$ on $Y$, a basic object of study in algebraic geometry is the graded ring
\begin{align}
R_{Y,D} &:=\dsum_{d \geq 0} H^0(Y,\sheaf_Y(dD)), \label{non-positive-sum}
\end{align}
and in some `natural' scenarios it turns out that $Y \cong \proj R_{Y,D}$ (e.g.\ see assertion \ref{wt-proj-compactn-to-filtrn} of Proposition \ref{prop-spec-completion}). Now let $X$ be the complement of the support of $D$ in $Y$ and $F_{D,d} := H^0(Y,\sheaf_Y(dD))|_X$. Then $\scrF_D := \{F_{D,d}\}_{d \in \zz}$ defines a filtration on $A := \Gamma(X, \sheaf_X)$, and it follows from our definition that $\af \cong R_{Y,D}$. Our choice of taking the sum over components with only non-negative components of $\scrF$ (in the definition of $\af$) comes from emulating the definition of $R_{Y,D}$. Initial supporting evidences for this choice come from Example \ref{general-weighted-example} and Proposition \ref{prop-spec-completion} below. 
\end{rem}

\begin{example}\label{general-weighted-example}
Let $\omega_1, \ldots, \omega_n \in \zz$ such that $\gcd(\omega_1, \ldots, \omega_n) = 1$ and
\begin{align}
\text{at least one of the $\omega_i$'s is positive.} \label{positive-condition}
\end{align}
Let $\omega$ be the {\em weighted degree} on $A := \kk[x_1, \ldots, x_n]$ corresponding to weights $\omega_k$ for $x_k$, $1 \leq k \leq n$. Define a filtration $\scrF_\omega$ on $A$ as $\scrF_\omega := \{f \in A: \omega(f) \leq d\}_{d \geq 0}$ and let $\bar X^{\scrF_\omega} := \proj \af$. If each $\omega_k > 0$, then we noted in the beginning of Section \ref{sec-intraborate} that $\bar X^{\scrF_\omega}$ is the classical weighted projective space $\pp^n(\kk;1,d_1,\ldots,d_n)$; in particular, $\bar X^{\scrF_\omega}$ is a projective completion of $X := \spec A \cong \ank$. But it is straightforward to check that even if some of the $\omega_k$'s are non-positive, as long as \eqref{positive-condition} is satisfied, $\bar X^{\scrF_\omega}$ is a quasi-projective variety containing $X$ as a dense open subset, $X_\infty := \bar X^{\scrF_\omega} \setminus X$ is an irreducible $n-1$ dimensional subvariety of $\bar X^{\scrF_\omega}$ and for every polynomial $f \in A$, $\omega(f)$ is precisely the {\em pole} of $f$ along $X_\infty$ (in the terminology of Section \ref{sec-semiquasintro}, $\omega$ is a {\em semidegree}, and the last assertion of the preceding sentence is a precursor of Proposition \ref{pole-and-degree}). Note that $\bar X^{\scrF_\omega}$ is a {\em generalized weighted projective space} in the sense of Definition \ref{general-weighted-defn}.
\end{example}

\begin{rem}
Note that \eqref{positive-condition} is a natural condition in our context, since the underlying goal in Example \ref{general-weighted-example} is to construct varieties that {\em contain} $\ank$, and $\omega$ is a {\em degree-like function} (Definition \ref{degree-like-defn}) which is `expected' to measure the {\em pole} of polynomials along the boundary to be constructed. For the resulting space to be `reasonable,' at least some of the polynomials should grow near infinity, i.e.\ they should have positive poles.  
\end{rem}

\begin{defn}[Generalized weighted projective spaces] \label{general-weighted-defn}
Pick $\omega_0, \ldots, \omega_n \in \zz$ satisfying \eqref{positive-condition}. We give two equivalent definitions for the {\em generalized weighted projective space} $\pp^n(\kk;\omega_0, \ldots, \omega_n)$ (note that our definitions agree with the definition of weighted projective spaces in the classically treated case that $\omega_j$'s are all positive). 
\end{defn}
{
\renewcommand{\thethm}{\ref{general-weighted-defn}.1}
\begin{defn}
Let $\omega$ be the weighted degree on $S := \kk[x_0, \ldots, x_n]$ corresponding to weights $\omega_k$ for $x_k$, $0 \leq k \leq n$. Then
\begin{align*}
&\pp^n(\kk;\omega_0, \ldots, \omega_n) := \proj \left(\sum_{d \geq 0} S_d \right),\quad \text{where}\\
&S_d := \{f \in S: f\ \text{weighted homogeneous with respect to}\ \omega,\ \omega(f) = d\}.
\end{align*}
\end{defn}

\renewcommand{\thethm}{\ref{general-weighted-defn}.2}
\begin{defn}
W.l.o.g.\ assume $\omega_0 > 0$. Pick linearly independent vectors $v_1, \ldots, v_n \in \rr^n$. Define $v_0 := -\frac{1}{\omega_0}\sum_{k=1}^n \omega_k v_k$ and set $L$ to be the $\zz$-span of $v_0, \ldots, v_n$ (i.e.\ $L$ is the {\em lattice} generated by $v_0, \ldots, v_n$). Then $\pp^n(\kk;\omega_0, \ldots, \omega_n)$ is the {\em toric variety} corresponding to the fan `determined by' $v_0, \ldots, v_n$ in the lattice $L$. More precisely, for each $i$, $0 \leq i \leq n$, let $\sigma_i$ be the ($n$-dimensional simplicial) cone generated by all $v_j$'s other than $v_i$, and let $\Sigma_{\omega_0, \ldots, \omega_n}$ be the simplicial fan consisting of those $\sigma_i$ (and all their sub-cones) for which $\omega_i > 0$. Then $\pp^n(\kk;\omega_0, \ldots, \omega_n)$ is the toric variety corresponding to the fan $\Sigma_{\omega_0, \ldots, \omega_n}$ and the lattice $L$. 

\xdefinecolor{c1}{HTML}{D4ACEE} 
\xdefinecolor{c2}{HTML}{EEc6AC} 
\xdefinecolor{c3}{HTML}{EDEEAC} 
\xdefinecolor{c4}{HTML}{B3EEAC}

\begin{figure}[h]
\begin{center}
\begin{tikzpicture}[scale=0.5]
	\begin{scope}[shift={(0,0)}]
		\draw [gray,  line width=0pt] (-3.5,-2.5) grid (2.5,2.5);	
		
		\foreach \xss in {.1,.2,...,2}{
 			\draw[color=c1] (\xss,0) -- (0,\xss); 
			\draw[color=c2] (\xss,0) -- (-1.5*\xss, -\xss); 
 			\draw[color=c3] (0,\xss) -- (-1.5*\xss, -\xss); 
 		}
 		\draw [blue, thick] (0,2) |- (2,0);
 		\draw [blue, thick, <->] (0,1) |- (1,0);
		\draw [blue] (1,0) node  [below] {$v_1$};
		\draw [blue] (0,1) node  [left] {$v_2$};
		\draw [blue,thick, ->] (0,0) -- (-3,-2) node [below] {$v_0$};
       	
		\draw (-.5,-3.5) node {\small{$\Sigma_{1,3,2}$}};    	

	\end{scope}
	
	\begin{scope}[shift={(7,0)}]
		\draw [gray,  line width=0pt] (-3.5,-2.5) grid (2.5,2.5);	
		
		\foreach \xss in {.1,.2,...,2}{
 			\draw[color=c1] (\xss,0) -- (0,\xss); 
 			\draw[color=c3] (0,\xss) -- (-1.5*\xss, \xss); 
 		}
 		\draw [blue, thick] (0,2) |- (2,0);
 		\draw [blue, thick, <->] (0,1) |- (1,0);
		\draw [blue] (1,0) node  [below] {$v_1$};
		\draw [blue] (0,1) node  [left] {$v_2$};
		\draw [blue,thick, ->] (0,0) -- (-3,2) node [left] {$v_0$};
       	
		\draw (-.5,-3.5) node {\small{$\Sigma_{1,3,-2}$}};  
	\end{scope}
	
	\begin{scope}[shift={(13,0)}]
		\draw [gray,  line width=0pt] (-2.5,-2.5) grid (3.5,2.5);	
		
		\foreach \xss in {.1,.2,...,3}{
 			\draw[color=c1] (\xss,0) -- (0,2*\xss/3); 
			\draw[color=c2] (\xss,0) -- (\xss, -2*\xss/3); 
 		}
 		\draw [blue, thick] (0,2) |- (3,0);
 		\draw [blue, thick, <->] (0,1) |- (1,0);	
		\draw [blue] (1,0) node  [below] {$v_1$};
		\draw [blue] (0,1) node  [left] {$v_2$};
		\draw [blue,thick, ->] (0,0) -- (3,-2) node [below] {$v_0$};
       	
		\draw (.5,-3.5) node {\small{$\Sigma_{1,-3,2}$}};  
	\end{scope}	
	
	\begin{scope}[shift={(20,0)}]
		\draw [gray,  line width=0pt] (-2.5,-2.5) grid (3.5,2.5);	
		
		\foreach \xss in {.1,.2,...,3}{
 			\draw[color=c1] (\xss,0) -- (0,2*\xss/3); 
 		}
 		\draw [blue, thick] (0,2) |- (3,0);
 		\draw [blue, thick, <->] (0,1) |- (1,0);
		\draw [blue] (1,0) node  [below] {$v_1$};
		\draw [blue] (0,1) node  [left] {$v_2$};
       	
		\draw (.5,-3.5) node {\small{$\Sigma_{1,-3,-2}$}};  
	\end{scope}		
	
\end{tikzpicture}
\caption{Fans of generalized weighted projective surfaces}  \label{fan-picture}
\end{center}
\end{figure}
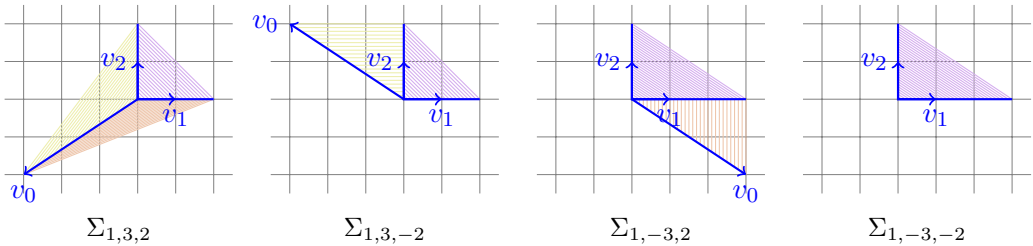
\end{defn}
\addtocounter{thm}{-2}
}

The following proposition establishes the basic dictionary between completions of an affine variety and filtrations on its coordinate ring.

\begin{prop}\label{prop-spec-completion}
Let $X$ be a (possibly reducible) affine variety. 
\begin{enumerate}
\item \label{filtrn-to-compactn} If $\scrF := \{F_d : d \in  \zz \}$ is a filtration on $\kk[X]$, then $\xf := \proj \kk[X]^{\scrF}$ is a reduced scheme containing $X$ as a (Zariski) open dense subset. The embedding $X \into \xf$ is induced via the map:
$$ X = \spec \kk[X] \ni \ppp \mapsto \dsum_{d \geq 0} (\ppp \cap F_d) \in \proj \kk[X]^\scrF = \xf.$$
The complement of $X$ in $\xf$ is the zero set of $(1)_1$ and is isomorphic (as a scheme) to $\proj \gring$. Moreover,
\begin{compactenum}
\item \label{non-trivially-finite-assumption} If $\scrF$ is finitely generated, then 
\begin{compactenum}
\item $\xf$ is a closed subvariety of a generalized weighted projective space, and
\item $\xf\setminus X$ is the support of an effective ample divisor on $\xf$.
\end{compactenum}
\item \label{projective-assumption} If in addition (to the assumptions of \ref{non-trivially-finite-assumption}) $F_0 = \kk$, then $\xf$ is a projective variety.
\end{compactenum}
\item  \label{compactn-to-filtrn}
\begin{enumerate}
\item \label{graded-compactn-to-filtrn} If $S$ is a graded $\kk$-algebra such that $\proj S$ contains $X$ as a dense open subset and $\proj S \setminus X = V(h)$ for some homogeneous $h \in X$ of degree $d > 0$, then there is a filtration $\scrF$ on $\kk[X]$ such that $\xf \cong \proj S$. 
\item \label{wt-proj-compactn-to-filtrn} If $Y$ is a closed subvariety of a generalized weighted variety containing $X$ such that $Y \setminus X$ is the support of an effective ample divisor, then there is a finitely generated filtration $\scrF$ on $\kk[X]$ such that $\xf \cong Y$.
\end{enumerate} 
\end{enumerate}
\end{prop}

\begin{proof}
We omit the proof of assertion \ref{filtrn-to-compactn}, since it follows in a straightforward matter from the definitions. For assertion \ref{graded-compactn-to-filtrn}, recall that $\proj S \setminus V(h) \cong \spec S_{(h)}$, where $S_{(h)}$ is the subring of the localization $S_h$ of $S$ consisting of homogeneous elements of degree zero. By assumption there is an isomorphism $\phi :S_{(h)} \cong \kk[X]$. For each $k \geq 0$, define
\begin{align*}
F_k := \phi(S_{kd}/h^{kd}) = \{g \in \kk[X] : \phi^{-1}(g) \in S_{kd}/h^{kd}\}.
\end{align*}
Then it is straightforward to check that the $\scrF := \{F_i\}_{i \geq 0}$ is a filtration on $\kk[X]$ which satisfies the claim of assertion \ref{graded-compactn-to-filtrn}. \\

Now we prove Assertion \ref{wt-proj-compactn-to-filtrn}. Since $Y$ is a closed subvariety of a generalized of a generalized weighted variety, it turns out that $Y = \proj S$ for a finitely generated $\nn$-graded $\kk$-algebra $S$ such that $S_0 = H^0(Y, \sheaf_Y)$ (this follows from arguments similar to the proof of the corresponding statements for subvarieties of projective spaces \cite[Proposition 5.15, Corollary II.5.16]{hart}). By assumption there is an effective Cartier divisor $D$ on $Y$ with $\supp D = Y\setminus X$ and global sections $f_0, \ldots, f_k$ of $\sheaf_Y(D)$ such that $\phi : Y \into \pp^k$ which maps $y \mapsto [f_0(y): \cdots : f_k(y)]$ induces an immersion. Since $D$ is effective, we may (and will) w.l.o.g.\ assume $f_0 = 1$. It follows that $\kk[X] = \kk[f_1, \ldots, f_k]$. Let $t_1, \ldots, t_k$ be indeterminates and define a filtration $\scrF = \{F_d\}_{d \geq 0}$ on $\kk[X]$ with  
\begin{align*}
F_d &:= \{f \in \kk[X]: f = G(f_1, \ldots, f_k)\ \text{for some}\ G \in S_0[t_1, \ldots, t_k],\ \deg(G) \leq d\},\quad d \geq 0.
\end{align*}
Note that $\scrF$ is generated by construction, since $\kk[X]^\scrF = S_0[(f_0)_1, \ldots, (f_k)_1]$. We now show that $\xf \cong Y$.\\

Let $U := \spec S_0$. Pick $\kk$-algebra generators $g_1, \ldots, g_l$ of $R$. Then $\xf$ is the closure in $U \times \pp^k$ of the image of the map $\tilde \phi: X \to U \times \pp^k$ given by $x \mapsto (g(x), \phi(x))$, where $g := (g_1, \ldots, g_l)$. Now note that $\tilde \phi$ extends naturally to all of $Y$ (since each $g_i \in H^0(Y, \sheaf_Y)$), and moreover, since $\phi$ is an immersion of $Y$ into $\pp^k$, it follows that $\tilde \phi$ induces an immersion of $Y$ into $\xf$. Now note that $S_0$ is the zeroth component of both $S$ and $\kk[X]^\scrF$, so that both $Y$ and $\xf$ are proper over the affine variety $U$. Since the morphism $\xf \to U$ is simply the extension of the morphism $Y \to U$, the arguments of \cite[Proposition 1]{goodman-affine} show that $\tilde \phi(Y) = \xf$, as required.  
\end{proof}

\begin{defn} \label{projinition}
We say that a filtration $\scrF$ on the coordinate ring of an affine variety $X$ is {\em projective} if it satisfies the hypothesis of assertion \ref{projective-assumption} of Proposition \ref{prop-spec-completion}. 
\end{defn}

Example \ref{general-weighted-example} is our most basic example of filtrations and corresponding immersions into weighted projective spaces. We now give some other examples.

\begin{example}
Let $X := \ank$, $n \geq 2$. Pick an integer $k$, $1 \leq k < n$, and define a filtration $\scrF := \{F_d\}_{d \geq 0}$ on $\kk[X] = \kk[x_1, \ldots, x_n]$ as follows: take $F_1$ to be the $\kk$-linear span of all monomials of degree less than or equal to two excluding those of the form $x_ix_j$ with $i \geq j > k$. Let $F_d = (F_1)^d$ for $d \geq 1$. Then $\xf$ is isomorphic to the blow-up of $\pp^n(\kk)$ along the subspace $V := V(x_0, \ldots, x_k)$. 
\end{example}

\begin{example}
If $X$ is a normal affine variety with trivial Picard group and $\bar{X} \subseteq \pNk$ be a normal projective completion of $X$ such that $X_\infty := \bar{X}\setminus X$ is irreducible, then $X_\infty$ is the support of an effective ample divisor and therefore $\bar X$ is isomorphic to $\xf$ for some filtration $\scrF$ on $\kk[X]$ (Proposition \ref{prop-spec-completion}).
\end{example} 

\begin{example} \label{curve-example}
If $\dim X = 1$, then every projective completion $Y$ of $X$ is isomorphic to $\xf$ for some filtration $\scrF$ on $\kk[X]$. Indeed, it follows from the Nakai-Moishezon criterion of ampleness (or the Riemann-Roch theorem) that $Y \setminus X$ is the support of an effective ample Cartier divisor \cite[Proposition 5]{goodman-affine}, so that the claim follows from Proposition \ref{prop-spec-completion}.
\end{example}

\begin{example} \label{surface-examples}
Every non-singular projective completion of affine surfaces comes from a filtration. More generally, it follows from \cite[Theorem 1]{goodman-affine} and the Nakai-Moishezon criterion that if $X$ is an arbitrary affine surface and $\bar X$ is a complete surface containing $X$ such that $\bar X \setminus X$ is the support of an effective Cartier divisor, then in fact $\bar X \setminus X$ is the support of an effective {\em ample} divisor, so that $\bar X$ is determined by a filtration on $\kk[X]$. Note that there exist complete normal surfaces $\bar X$ and open affine subsets $X$ of $\bar X$ such that $\bar X \setminus X$ is not the support of any Cartier divisor (see e.g.\ \cite{overflowing-non-amply}) - so that $\bar X$ is {\em not} determined by any filtration on $\kk[X]$.
\end{example}

\begin{example}
There are {\em non-singular} projective completions of (non-singular) affine varieties of dimension $\geq 3$ which are {\em not} determined by a filtration, as the following example of Hironaka shows. Let $\phi:V' \to V$ be a morphism of 3-dimensional projective varieties, with $V'$ non-singular and $V$ non-singular except at a point $P$, such that $\phi^{-1}(P)$ is a curve and $\phi$ is an isomorphism on the rest of $V'$ and $V$. Let $W$ be a hyperplane section of $V$ through $P$. Then $X := V' \setminus \phi^{-1}(W) \cong V \setminus W$ is affine, but $\phi^{-1}(W)$ is {\em not} the support of any effective ample divisor on $V'$ (see \cite[Section 1]{goodman-affine}). Therefore Proposition \ref{prop-spec-completion} implies that there is no filtration $\scrF$ on $A$ such that $\xf \cong V'$.
\end{example}

\section{\Gf and \Sgf} \label{sec-semiquasintro}
\newcommand{\adelta}{A^\delta}
\newcommand{\kxdelta}{\kk[X]^\delta}

\begin{defn} \label{degree-like-defn}
A {\em degree-like function} on $A$ is a map $\delta: A\setminus \{0\} \to \zz \cup \{-\infty\}$ such that:
\begin{compactenum}
\item \label{degreelike-prop-0} $\delta(\kk) = 0$.
\item \label{degreelike-prop-1} $\delta(f+g) \leq \max\{\delta(f), \delta(g)\}$ for all $f, g \in A$, with $<$ in the preceding equation implying $\delta(f) = \delta(g)$.
\item \label{degreelike-prop-2} $\delta(fg) \leq \delta(f) + \delta(g)$ for all $f, g \in A$.
\end{compactenum}
\end{defn}

\begin{rem}
We allow degree-like functions to have $-\infty$ as a value, since they do correspond to some filtrations (via the correspondence defined below). On the other hand, we will see (from Proposition \ref{finite-prop} and Theorem \ref{thm-noeth-integral-subdegree} below) that the degree-like functions associated to completions of affine varieties are always integer-valued.
\end{rem}

There is a one-to-one correspondence between degree-like functions and filtrations:
\begin{center}
\begin{tabular}{ccc}
{\em Filtrations} & $\longleftrightarrow$ & {\em degree-like functions} \\
$\scrF	= \{F_d\}_{d \in \zz}$											& $\longrightarrow$	& $\delta_\scrF:f \in A \mapsto \inf\{d:f \in F_d\}$ \\
$\scrF_\delta := \{F_d := \{f \in A: \delta(f) \leq d\}\}_{d \in \zz}$ 	& $\longleftarrow$ 	& $\delta$
\end{tabular}
\end{center}

In the remainder of this article we identify degree-like functions with the corresponding filtrations. In particular, we refer to a degree-like function $\delta$ as {\em projective} (resp.\ {\em finitely generated}) iff the corresponding filtration $\scrF_\delta$ is projective (resp.\ finitely generated). Similarly, $\kxdelta$, $\gr \kxdelta$ and $\xdelta$ will be shorthand notations respectively for $\kk[X]^{\scrF_\delta}$, $\gr \kk[X]^{\scrF_\delta}$ and $\bar X^{\scrF_\delta}$. \\

The protagonists of this article are two classes of degree-like functions which satisfy stronger versions of the multiplicative property (i.e.\ property \ref{degreelike-prop-2} above).

\begin{defn}
\mbox{}
\begin{compactitem}
\item A degree-like function $\delta$ on $A$ is a {\em \gengf} iff $\delta(fg) = \delta(f) + \delta(g)$ for all $f,g \in A\setminus\{0\}$.
\item We say that $\delta$ is a {\em \gensgff} if there are \gengff s $\delta_1, \ldots, \delta_N$ such that 
\begin{align} \label{gensgf-condition}
	\delta(f) = \max_{1 \leq i \leq N} \delta_{i}(f) \quad \text{for all}\ f \in A\setminus\{0\}.
\end{align}
\end{compactitem}

Given a \gensgf $\delta$ as in \eqref{gensgf-condition}, we may assume by getting rid of some $\delta_i$'s if necessary that every $\delta_i$ that appears in \eqref{gensgf-condition} is {\em not redundant}, i.e.\ for every $i$, there is an $f \in A$ such that $\delta_{i}(f) > \delta_{j}(f)$ for all $j \neq i$. In this case we will say that \eqref{gensgf-condition} is a {\em minimal presentation} of $\delta$.  
\end{defn}

\begin{example}[Weighted Degree] \label{weighted-semi-example}
As noted in the introduction, every weighted degree $\delta$ on the polynomial ring $\kk[x_1, \ldots, x_n]$ is a semidegree. 
\end{example}

\begin{example}[Subdegrees determined by rational polytopes]\label{toric-example}
Let $X$ be the $n$-torus $(\kk^*)^n$ and $A := \kk[x_1, x_1^{-1}, \ldots, x_n, x_n^{-1}]$ be its coordinate ring. Let $\scrP$ be a convex rational polytope (i.e.\ a convex polytope in $\rr^n$ with vertices in $\qq^n$) of dimension $n$ containing origin in its interior. Recall (from introduction) that $\scrP$ determines a degree-like function on $A$. More precisely, recall that if $\delta_\scrP: A \setminus\{0\} \to \qq$ is defined as
$$\delta_\scrP(\sum a_\alpha x^\alpha) := \inf\{r \in \rr: r \geq 0,\ \alpha \in r\scrP\ \text{for all}\ \alpha \in \zz^n\ \text{such that}\ a_\alpha \neq 0\},$$ 
and $e$ is any positive integer such that $e\delta_\scrP$ is integer-valued, then $e\delta_\scrP$ is a degree-like function on $A$. We claim that $e\delta_\scrP$ is in fact a subdegree on $A$. Indeed, for each facet $\scrQ$ of $\scrP$, let $\omega_\scrQ$ be the smallest `outward pointing' integral vector normal to $\scrQ$ and let $c_\scrQ := \langle \omega_\scrQ, \alpha_\scrQ \rangle$, where $\alpha_\scrQ$ is any element of the hyperplane that contains $\scrQ$. Define 
$$\delta_\scrQ(\sum a_\alpha x^\alpha) := \max_{a_\alpha \neq 0} \frac{\langle \omega_\scrQ, \alpha \rangle}{c_\scrQ}.$$
Then it is straightforward to see that $e\delta_\scrQ$ is a weighted degree on $A$ for each facet $\scrQ$ of $\scrP$, and $e\delta_\scrP = \max\{e\delta_\scrQ: \scrQ$ is a facet of $\scrP\}$, so that $e\delta_\scrP$ is a subdegree, as claimed. It is not hard to see that $X^{\delta_\scrP}$ is isomorphic to the {\em toric variety} $X_\scrP$ determined by $\scrP$.
\end{example}

\begin{example}[An iterated semidegree]\label{iterated-semi-baby-example}
Let $X := \affine{2}{\kk}$ and $A := \kk[x_1,x_2]$. Define a filtration $\scrF := \{F_d: d \geq 0\}$ on $A$ by setting $F_0 := \kk~$, $F_1 := \kk\langle 1, x_1^2-x_2^3 \rangle~$, $F_2 := (F_1)^2 + \kk\langle x_2 \rangle~$, $F_3 := F_1F_2 + \kk\langle x_1 \rangle$ and $F_d := \sum_{j=1}^{d-1} F_jF_{d-j}$ for $d \geq 4$. Let $\delta := \delta_\scrF$. Note that $\delta(x_1) = 3$ and $\delta(x_2)= 2$, but $\delta(x_1^2 - x_2^3) = 1$. We claim that $\delta$ is a semidegree. Indeed, 
\begin{align}
\adelta = \kk[(1)_1, (x_1)_3, (x_2)_2, (x_1^2-x_2^3)_1] \cong \kk[X_1,X_2,Y,Z]/\langle YZ^5 - X_1^2 +X_2^3 \rangle\: , \label{iterated-graded}
\end{align}
where the last isomorphism is induced by a $\kk$-algebra homomorphism which sends \mbox{$X_1 \mapsto (x_1)_3$}, $X_2 \mapsto (x_2)_2$, $Y \mapsto (x_1^2-x_2^3)_1$ and $Z \mapsto (1)_1$. The inverse image of the ideal $I := \langle (1)_1 \rangle$ of $\adelta$ under this isomorphism coincides with ideal $\langle Z, YZ^5 - X_1^2 +X_2^3 \rangle = \langle Z, X_1^2  - X_2^3 \rangle$. Since the latter is a prime ideal of $\kk[X_1,X_2,Y,Z]$, it follows that $I$ is a prime ideal of $\adelta$ as well. Theorem \ref{gengfsgf-characterization} then implies that $\delta$ is a semidegree. Identity \eqref{iterated-graded} implies that $\xdelta$ is isomorphic to the hypersurface $V(YZ^5 - X_1^2 +X_2^3)$ of the weighted projective space $\pp^3(\kk;1,3,2,1)$ with (weighted homogeneous) coordinates $[Z:X_1:X_2:Y]$. The semidegree $\delta$ is an example of an {\em iterated semidegree} (to be introduced in the forthcoming continuation of this article).
\end{example}

\section{Basic Structure of Subdegrees}\label{sec-sub-structure}
\renewcommand{\filtrationring}{\ensuremath{A}}
\renewcommand{\filtrationchar}{\ensuremath{\delta}}
\newcommand{\abarf}{\ensuremath{\profingg{A}{\bar \scrF}}}
\newcommand{\tildeAEta}{\ensuremath{\profingg{\tilde A}{\eta}}}
\newcommand{\tildeABarEta}{\ensuremath{\profingg{\tilde A}{\bar \eta}}}
\newcommand{\afd}{(\af)^{[d]}}
\newcommand{\del}[1]{(#1)_{\delta(#1)}}
\newcommand{\locali}{\sheaf_{V_{i}, \xdelta}}
\newcommand{\localj}{\sheaf_{V_{j}, \xdelta}}
The goal of this section is to prove the following theorem and its corollaries. Recall that an ideal $I$ of a ring $R$ is {\em decomposable} if it is the intersection of finitely many primary ideals of $R$. 

\begin{thm}[Structure theorem for subdegrees]\label{gengfsgf-characterization}
Let $\delta$ be a degree like function on an integral domain $A$ and $I$ be the ideal of $\adelta$ generated by $(1)_1$ (recall the notation from Section \ref{sec-filtrintro} that each for $f \in A$ and each $d \geq \delta(f)$, we denote the `copy' of $f$ in the $d$-th graded component of $\adelta$ by $(f)_d$).
\begin{compactenum}
\item\label{gengf-characterization} $\delta$ is a \gengf if and only if $I$ is a prime ideal.
\item\label{gensgf-characterization} $\delta$ is a \gensgf if and only if $I$ is a decomposable radical ideal. 
\item \label{sgf-unique-decomp} If $\delta$ is a \gensgff, then \gengff s $\delta_1, \ldots, \delta_N$ of a minimal presentation $\delta = \max_{1\leq i\leq N} \delta_i$ of $\delta$ are unique.
\item\label{noetherian-gensgf-characterization} If $\profing$ is Noetherian, then $\delta$ is a \gensgf if and only if $I$ is a radical ideal.
\end{compactenum}
\end{thm}

Theoem \ref{gengfsgf-characterization} was proven in my Ph.D. thesis \cite[Section 2.2]{pisis}. Later I came to know of the article \cite{szpiro} (I heartily thank Professor Bernard Teissier for bringing it to my attention) - and assertions \ref{gengf-characterization}, \ref{gensgf-characterization} and \ref{noetherian-gensgf-characterization} of Theorem \ref{gengfsgf-characterization} follow from the results of \cite{szpiro} (see Remark \ref{szpiro-sgf} below). However, I present here a complete proof of Theorem \ref{gengfsgf-characterization}, since proving only assertion \ref{sgf-unique-decomp} takes almost the same amount of work. The essential ingredients are Lemmas \ref{gensgf-characterize-lemma} and \ref{limit-semidegree}. 

\begin{cor} \label{sgf-reformulation}
Assume $\adelta$ is Noetherian. Then $\delta$ is a \gensgf on $A$ iff $\delta(f^k) = k\delta(f)$ for all $f \in A$ and $k \geq 0$.
\end{cor}

\begin{proof}
Let $I$ be the ideal of $\adelta$ from the statement of Theorem \ref{gengfsgf-characterization}. Note that for every $f \in A$ and $d \in \zz$, $(f)_d \in I$ iff $\delta(f) < d$. It follows that $I$ is radical iff $\delta(f^k) = k\delta(f)$ for all $f \in A$ and $k \geq 0$. Now the corollary follows from assertion \ref{noetherian-gensgf-characterization} of Theorem \ref{gengfsgf-characterization}.
\end{proof}

\begin{rem} \label{szpiro-sgf}
\cite{szpiro} considers {\em order functions}, which are {\em negative} of our degree-like functions with the exception that the values are allowed to be in $\rr \cup \infty$ (instead of $\zz \cup \infty$). An order function $\nu$ on a ring $A$ is {\em homogeneous} if $\nu(f^k) = k\delta(f)$ for all $k \geq 0$ and $f \in A$. \cite[Theorem 1]{szpiro} states that $\nu$ is homogeneous iff $\nu$ is the maximum of a (possibly infinite) collection of {\em valuations}. Assertions \ref{gengf-characterization}, \ref{gensgf-characterization} and \ref{noetherian-gensgf-characterization} of Theorem \ref{gengfsgf-characterization} and Corollary \ref{sgf-reformulation} follow from this result. However, the results of \cite{szpiro} (or arguments of their proofs) does not seem to imply Assertion \ref{sgf-unique-decomp} of Theorem \ref{gengfsgf-characterization}.
\end{rem}

\begin{unorderedproof}[Proof of Theorem \ref{gengfsgf-characterization}]{thm}{1}
Assertion \ref{gengf-characterization} follows from the definition of semidegrees and the observation that for all $f \in A$ and $d \in \zz$, $\delta(f) = d$ iff $(f)_d \in \adelta \setminus I$ . Assertion \ref{noetherian-gensgf-characterization} is a consequence of assertion \ref{gensgf-characterization}. Therefore it remains to prove assertions \ref{gensgf-characterization} and \ref{sgf-unique-decomp}. We now prove the `only if' direction of assertion \ref{gensgf-characterization}.\\

Assume $\delta$ is a subdegree. Let $\delta = \max_{1\leq i\leq N} \delta_i$ be a minimal presentation of $\delta$. For each $i$ with $0 \leq i \leq N$, fix an $f_i \in A$ such that 
\begin{align}
d_i  &:= \delta_i(f_i) > \delta_j(f_i)\ \text{for}\ 1\leq j \neq i \leq N. \label{f_i-property}
\end{align}
Then, in particular, $\delta(f_i) = d_i \in \zz$. Assertion \ref{part2-gensgf-characterize-lemma} of the following lemma implies the `only if' direction of assertion \ref{gensgf-characterization} of the theorem.

\begin{prolemma}\label{gensgf-characterize-lemma}
For each $i$, $1 \leq i \leq N$, let $\ppp_i := (I:(f_i)_{d_i}) = \{F \in \adelta : F(f_i)_{d_i} \in I\}$.
\begin{compactenum}[(1)]
\item \label{part1-gensgf-characterize-lemma} The collection of homogeneous elements of $\ppp_i$ is $L_i := \{(f)_{d}: f \in A,\ \delta_i(f) < \delta(f) \leq d\}.$
\item \label{part2-gensgf-characterize-lemma} For each $i$, $\ppp_i$ is a distinct minimal prime ideal of $\profing$ containing $I$. In particular, $(f_i)_{d_i} \in \left(\bigcap_{j \neq i}\ppp_j \right) \setminus \ppp_i$ for each $i$. Moreover, $I = \bigcap_{i=1}^N \ppp_i$.
\end{compactenum}
\end{prolemma}

\begin{proof}
We first prove assertion \ref{part1-gensgf-characterize-lemma}. Fix an $i$, $1 \leq i \leq N$. Let $(f)_d$ be an arbitrary homogeneous element of $\ppp_i$. We will show that $(f)_d \in L_i$. If $\delta(f) < d$, then $\delta_i(f) \leq \delta(f) < d$ and $(f)_d \in L_i$. So assume $\delta(f) = d$. Since $(f)_d(f_i)_{d_i} = (ff_i)_{d+d_i} \in I$, it follows that $\delta(ff_i) < d + d_i$. But then $\delta_i(ff_i) = \delta_i(f)+\delta_i(f_i) < d + d_i$, so that $\delta_i(f) < d = \delta(f)$, and thus $(f)_d \in L_i$. To summarize, all homogeneous elements of $\ppp_i$ belong to $L_i$. \\

Now we show that $L_i \subseteq \ppp_i$. Pick $(f)_d \in L_i$. If $d > \delta(f)$, then $(f)_d \in I \subseteq \ppp_i$. So assume $d = \delta(f)$. Then $\delta_i(f) < d = \delta(f)$, and thus $\delta_i(ff_i) = \delta_i(f) + \delta_i(f_i) < d + d_i$. Also, for each $j \neq i$, $\delta_j(f_i) < d_i$, so that $\delta_j(ff_i) = \delta_j(f) + \delta_j(f_i) < d + d_i$. It follows that $\delta(ff_i) < d + d_i$, so that $(f)_d(f_i)_{d_i} = (ff_i)_{d+d_i} \in I$. Therefore $(f)_d \in \ppp_i$, which proves that $L_i \subseteq \ppp_i$ and completes the proof of assertion \ref{part1-gensgf-characterize-lemma} of the lemma.\\

To see that $\ppp_i$ is prime, let $(g_1)_{e_1}, (g_2)_{e_2}$ be homogeneous elements of $\adelta$ such that $(g_1)_{e_1}(g_2)_{e_2} \in \ppp_i$. Since $\ppp_i$ is a homogeneous ideal, it suffices to show that one of the $(g_j)_{e_j}$'s belong to $\ppp_i$. Due to assertion \ref{part1-gensgf-characterize-lemma} of the lemma, $(g_1g_2)_{e_1 + e_2} = (g_1)_{e_1}(g_2)_{e_2} \in L_i$, which means that $\delta_i(g_1g_2) = \delta_i(g_1) + \delta_i(g_2) < e_1 + e_2$. Since $e_j \geq \delta(g_j) \geq \delta_i(g_j)$ for each $j$, it follows that there is at least one $j$ such that $\delta_i(g_j) < e_j$. Then assertion \ref{part1-gensgf-characterize-lemma} of the lemma implies that $(g_j)_{e_j} \in \ppp_i$, as required.\\

Next we show that $\ppp_i$'s are distinct minimal prime ideals of $\adelta$ containing $I$. Indeed, that $\ppp_i$'s are distinct follows from assertion \ref{part1-gensgf-characterize-lemma} and the observation that $(f_i)_{d_i} \in (\bigcap_{j \neq i}L_j) \setminus L_i$ for $1 \leq i \leq N$. Now fix $i$, $1 \leq i \leq N$, and a prime ideal $\ppp$ of $\profing$ such that $I \subseteq \ppp \subseteq \ppp_i$. Then $(f_i)_{d_i} \not\in \ppp$ (since assertion \ref{part1-gensgf-characterize-lemma} implies $(f_i)_{d_i} \not\in \ppp_i$). But if $(g)_e \in \ppp_i$, then $(g)_e(f_i)_{d_i} \in I \subseteq \ppp$ and it follows that $(g)_e \in \ppp$. Hence $\ppp_i \subseteq \ppp$, and therefore $\ppp_i = \ppp$. Consequently, $\ppp_i$ is a minimal prime ideal of $\adelta$ containing $I$, as required.\\

Finally, pick any homogeneous $(g)_e \in \bigcap_{i=1}^N \ppp_i$. Then by assertion \ref{part1-gensgf-characterize-lemma}, for each $i$, $\delta_i(g) < e$. Therefore $\delta(g) = \max_i \delta_i(g)< e$, and hence $(g)_e \in I$. It follows that $\bigcap_{i=1}^N \ppp_i = I$, which concludes the proof of the lemma.
\end{proof}

Now we prove assertion \ref{sgf-unique-decomp} of the theorem (the `if' direction of assertion \ref{gensgf-characterization} will be proven at the end). Let $\delta = \max_{1\leq i\leq N'} \delta'_i$ be another minimal presentation of $\delta$. Then there exist $f'_1, \ldots, f'_{N'} \in A$  such that
\begin{align}
d'_i &:= \delta'_i(f'_i) > \delta'_j(f'_i)\ \textrm{for}\ 1\leq j \neq i \leq N'. \tag{\oldref{f_i-property}$'$} \label{f'_i-property}
\end{align}
Then $d'_i$ are integers. Let $\ppp'_i := (I: (f'_i)_{d'_i})$, $1 \leq i \leq N'$. Then Lemma \ref{gensgf-characterize-lemma} implies that both $I = \bigcap_{i=1}^{N} \ppp_i$ and $I = \bigcap_{i=1}^{N'} \ppp'_i$ are minimal primary decompositions of $I$. The uniqueness of minimal primary decompositions then implies that $N=N'$ and, after an appropriate re-indexing of $\delta'_i$'s if necessary, $\ppp_i = \ppp'_i$ for all $i$, $1 \leq i \leq N$. \\

Fix an $i$, $1 \leq i \leq N$. Lemma \ref{gensgf-characterize-lemma} implies that 
$$(f'_i)_{d'_i} \in \bigcap_{j \neq i} L'_j \setminus L'_i = \bigcap_{j \neq i} L_j \setminus L_i\: ,$$ 
where $L_j$'s (resp.\ $L'_j$'s) are the set of homogeneous elements of $\ppp_j$'s (resp.\ $\ppp'_j$'s). 
It follows from the description of $L_j$'s and $L'_j$'s from assertion \ref{part1-gensgf-characterize-lemma} of Lemma \ref{gensgf-characterize-lemma} that for each $i$, $1 \leq i \leq N$,
\begin{align}
d'_i = \delta'_i(f'_i) = \delta(f'_i) = \delta_i(f'_i),\ \text{and} \label{delta-i-f'-i} \\
\delta_j(f'_i) < \delta_i(f'_i)\ \text{for all}\ j \neq i. \label{delta-j-f'-i}
\end{align}
Identity \eqref{delta-i-f'-i} and inequality \eqref{delta-j-f'-i} imply that property \eqref{f_i-property} remains true even after replacing $f_i$'s and $d_i$'s by $f'_i$'s and $d'_i$'s. Since $f_i$'s were assumed to be arbitrary elements in $A$ such that \eqref{f_i-property} is true, we may assume without loss of generality that $f_i = f'_i$ for each $i$. But then the following lemma implies that $\delta_i = \delta'_i$ for each $i$, as required.

\begin{prolemma} \label{limit-semidegree}
If $\delta = \max_{1\leq i\leq N} \delta_i$ is a minimal presentation for $\delta$ and $f_1, \ldots, f_N \in A$ satisfy \eqref{f_i-property}, then 
\begin{align}
\delta_i(f) = \lim_{k \to \infty} \delta((f_i)^kf) - \delta((f_i)^k) \label{defn-delta_i-from-f_i}
\end{align}
for all $f \in A$ and all $i$, $1 \leq i \leq N$.
\end{prolemma}

\begin{proof}
Fix an $i$, $1 \leq i \leq N$, and define $\tilde \delta_i(f) := \lim_{k \to \infty} \delta((f_i)^kf) - \delta((f_i)^k)$ for all $f \in A$. We first show that $\tilde \delta_i$ is well defined. Indeed, Let $f \in A$ and $k \geq 1$. Since $(f_i)_{d_i} \not\in I$ and $I$ is radical (assertion \ref{part2-gensgf-characterize-lemma} of Lemma \ref{gensgf-characterize-lemma}), it follows that $((f_i)_{d_i})^k = ((f_i)^k)_{kd_i} \not\in I$, so that $\delta((f_i)^k) = kd_i$. Therefore,
\begin{align}
\begin{split} \label{delta-i-seq-decreasing}
\delta((f_i)^{k+1}f) - \delta((f_i)^{k+1}) &\leq \delta((f_i)^kf) + \delta(f_i) - \delta((f_i)^{k+1})\\
		     			   &= \delta((f_i)^kf) + d_i - (k+1)d_i \\
		     			   &= \delta((f_i)^kf) - kd_i \\
		     			   &= \delta((f_i)^kf) - \delta((f_i)^k).
\end{split}
\end{align}
It follows that $\tilde \delta_i(f)$ is a well defined element in $\zz \cup \{-\infty\}$. \\

We now show that $\tilde \delta_i = \delta_i$. Note that for all $k \geq 0$ and all $f \in A$, $\delta((f_i)^kf) - \delta((f_i)^k) \geq \delta_i((f_i)^kf) - \delta((f_i)^k) = \delta_i(f)$ (the last equality uses that $\delta(f_i) = \delta_i(f_i)$), so that $\tilde \delta_i \geq \delta_i$. To see the opposite inequality, let $f \in A$ and $d \in \zz$ be such that $d \geq \delta_i(f)$. Then $kd_i + d \geq \delta_i((f_i)^kf)$ for all $k$. Moreover, \eqref{f_i-property} implies that for sufficiently large $k$, $kd_i + d > k\delta_j(f_i) + \delta_j(f) = \delta_j((f_i)^kf)$ for all $j \neq i$. It follows that for sufficiently large $k$, $\delta((f_i)^kf) \leq kd_i + d$, and hence $\delta((f_i)^kf) - \delta((f_i)^k) \leq d$. Since $\{\delta((f_i)^kf) - \delta((f_i)^k)\}$ is a decreasing sequence due to \eqref{delta-i-seq-decreasing}, it follows that $\tilde \delta_i(f) \leq d$. Summarizing, we have proved that $\tilde \delta_i(f) \leq d$ whenever $\delta_i(f)\leq d$. Applying the preceding statement to $d = \delta_i(f)$ when $\delta_i(f) \in \zz$, and otherwise letting $d$ converge to $\delta_i(f) = -\infty$ from above, we see that $\tilde \delta_i(f) \leq \delta_i(f)$. It follows that $\delta_i(f) = \tilde \delta_i(f)$, which completes the proof of the lemma and therefore, assertion \ref{sgf-unique-decomp} of Theorem \ref{gengfsgf-characterization}.
\end{proof}

It remains to prove the `if' direction of assertion \ref{gensgf-characterization} of the theorem. So assume $I$ is a decomposable radical ideal of $\profing$. We have to show that $\delta$ is a subdegree. It follows from the general theory of primary decomposition (e.g.\ as in \cite[Chapter 4]{am}) that there exist $(f_1)_{d_1}, \ldots, (f_N)_{d_N} \in \adelta \setminus I$ such that $\ppp_i := (I:(f_i)_{d_i})$ is a prime ideal containing $I$ for each $i$ and $I = \bigcap_{i=1}^N \ppp_i$. Since $(f_i)_{d_i} \not\in I$, it follows that $\delta(f_i) = d_i \in \zz$ for each $i$. For each $i=1, \ldots, N$, define $\delta_i:A \to \zz \cup \{-\infty\}$ as in \eqref{defn-delta_i-from-f_i}, i.e.\ 
\begin{align*} 
\delta_i(f) &:= \lim_{k \to \infty} \delta((f_i)^kf) - \delta((f_i)^k)
\end{align*}
for all $f \in A$. The same arguments as in the first part of the proof of Lemma \ref{limit-semidegree} imply that the relations of \eqref{delta-i-seq-decreasing} remain true and consequently $\delta_i$ is well defined for all $i$, $1 \leq i \leq N$. It follows from elementary but somewhat long computations that $\delta_i$'s satisfy the properties of a semidegree - we refer to the proof of \cite[Theorem 2.2.1]{pisis} for details. We now show that $\delta = \max_{i=1}^N \delta_i$. Indeed, fix $f \in A$ and $1 \leq i \leq N$. Then $\delta(f) \geq \delta(f_if) - \delta(f_i) \geq \delta_i(f)$ (the first inequality follows from the multiplicative property of degree-like functions and the second one is a consequence of \eqref{delta-i-seq-decreasing}). It follows that $\delta \geq\max_{i=1}^N \delta_i$. \\

For the opposite inequality, pick any $f \in A$. If $\delta(f) = -\infty$, then for all $i$ and for all $k \geq 0$, $\delta((f_i)^kf) = -\infty$, and therefore $\delta_i(f) = -\infty = \delta(f)$. So assume $\delta(f) \in \zz$. Since $(f)_{\delta(f)} \not\in I$, there is an $i$ such that $(f)_{\delta(f)} \not\in \ppp_i$. Note that $(f_i)_{d_i}$ also does not lie in $\ppp_i$, for otherwise $((f_i)_{d_i})^2$ would be an element of $I$ and this would imply that $(f_i)_{d_i} \in I$ (since $I$ is a radical ideal), contradicting our choice of $(f_i)_{d_i}$. Thus neither of $(f)_{\delta(f)}$ and $(f_i)_{d_i}$ is an element of $\ppp_i$. Since $\ppp_i$ is prime, it follows that for all $k \geq 0$, $((f_i)_{d_i})^k(f)_{\delta(f)} = ((f_i)^kf)_{kd_i+\delta(f)} \not\in \ppp_i$. Consequently $\delta(f_i^kf) = \delta(f) + kd_i = \delta(f) + \delta(f_i^k)$ for all $k \geq 0$. It then follows that $\delta_i(f) = \delta(f)$. Combining with the inequality proved in the preceding paragraph, this shows that $\delta =\max_{i=1}^N \delta_i$ and therefore $\delta$ is indeed a subdegree. This completes the proof of the theorem.
\end{unorderedproof}

\begin{cor} \label{sgf-components-at-infinity}
Let $\delta$ be a non-negative subdegree on $A$ and $\delta = \max_{i=1}^N \delta_i$ be its minimal presentation. Let $X := \spec A$. Then the number of the irreducible components of $X_\infty:=\xdelta \setminus X$ is $N$ if none of the $\delta_i$'s is the {\em zero} degree-like function and $N-1$ otherwise.
\end{cor}

\begin{proof}
Proposition \ref{prop-spec-completion} implies that $X_\infty = V(I) \subseteq \proj \adelta$, where $I$ is the ideal in $\profing$ generated by $(1)_1$. According to Lemma \ref{gensgf-characterize-lemma}, $I$ has a minimal prime decomposition of the form $I = \bigcap_{i=1}^N \ppp_i$, with each $\ppp_i$ being a prime ideal corresponding to $\delta_i$, $1 \leq i \leq N$. Consequently $X_\infty = V(I) = \bigcup_{i=1}^N V(\ppp_i)$. The corollary now follows from the following

\begin{claim*}
For each $i$, $1 \leq i \leq N$, the semidegree $\delta_i$ is identically zero on $A$ iff $V(\ppp_i) = \emptyset$.
\end{claim*}

\begin{proof}
Let $\adelta_+ = \dsum_{d \geq 1} \{(f)_d: f \in A\}$  be the {\em irrelevant ideal} of $\adelta$. Fix $i$, $1 \leq i \leq N$. Then
\begin{align*}
V(\ppp_i) = \emptyset & \iff\ \ppp_i \supseteq \adelta_+\\
			 & \iff\ \{(f)_{d}: f \in A,\ d > \delta_i(f)\} \supseteq \{(f)_d: f\in A,\ d \geq 1\}\ \text{(according to Lemma \ref{gensgf-characterize-lemma})}\\
			 & \iff\ \delta_i(f) \leq 0\ \text{for all}\ f \in A \\
			 & \iff\ \delta = \max \left(\{\delta_j: j \neq i,\ 1 \leq j \leq N\} \cup \{\delta_0\} \right)\ \text{(since}\ \delta\ \text{is non-negative)},
\end{align*}
where $\delta_0$ is the zero degree-like function on $A$, i.e.\ $\delta_0(f) := 0$ for all $f \in A$. Since $\delta_0$ is a semidegree and the minimal decomposition of $\delta$ is unique, it follows that $V(\ppp_i) = \emptyset \iff \delta_i = \delta_0$, which completes the proof of the claim.			 
\end{proof}
\renewcommand{\qedsymbol}{}
\end{proof}

\begin{example} \label{toric-number-of-components-example}
Applied to the subdegree of Example \ref{toric-example} determined by a convex rational polytope $\scrP$, Corollary \ref{sgf-components-at-infinity} implies the (standard) fact that the irreducible components of the complement of the torus in the toric variety $X_\scrP$ are in a one-to-one correspondence with the facets of $\scrP$.
\end{example}

\section{Normality at infinity} \label{sec-sub-normality}
In this section using the theory of {\em Rees' valuations} we work out the precise relation between the semidegrees in the minimal presentation of a subdegree and the orders of vanishing along the components of the hypersurface at infinity on the corresponding completion. We also establish our `Main Existence Theorem' which states that given a completion $\bar X$ of an affine variety $X$ determined by an arbitrary degree-like function, there is a subdegree on $\kk[X]$ which determines the {\em normalization at infinity} of $\bar X$ with respect to $X$. As a consequence we prove a `finiteness' property of divisorial valuations. We start with a presentation (following \cite[chapter XI]{mcdivisors}) of the relevant results of Rees (starting with a reminder of the notion of a {\em Krull domain}).

\begin{defn*}
A domain $B$ is a {\em Krull domain} iff 
\begin{compactenum}
\item $B_\ppp$ is a discrete valuation ring for all height one prime ideals $\ppp$ of $B$, and
\item every non-zero principal ideal of $B$ is the intersection of a finite number of primary ideals of height one.
\end{compactenum}
\end{defn*}

Every normal Noetherian domain is a Krull domain \cite[Section 41]{matsulgebra}. In particular, the integral closure of $\adelta$ is a Krull domain provided that $\adelta$ is finitely generated.\\

For an ideal $I$ of a ring $R$ define $\nu_I: R \to \nn \cup \{\infty\}$ and $\bar \nu_I: R \to \qq_+ \cup \{\infty\}$ by: 
\begin{align}
\begin{split} \label{nu-nu-bar}
\nu_I(x) &:= \sup\{m: x \in I^m\},\ \text{and} \\
\bar \nu_I(x) &:= \lim_{m \to \infty} \frac{\nu_I(x^m)}{m},
\end{split}
\end{align}
for all $x \in R$. Recall that the {\em integral closure} $\bar J$ of an ideal $J$ of $R$ is the ideal defined by: $\bar J := \{x \in R: x$ satisfies an equation of the form: $x^s + j_1x^{s-1}+ \cdots + j_s = 0$ with $j_k \in J^k$ for all $k = 1, \ldots, s\}$.

\setcounter{thm}{-1}
\begin{thm}[Rees' Theorem - see {\cite[Propositions 11.1 -- 11.6]{mcdivisors}}] \label{rees-thm}
For any ring $R$ and any ideal $I$ of $R$, $\bar \nu_I$ is well defined. Assume $R$ is a Noetherian domain. Then
\begin{asparaenum}[(1)]
\item \label{normalization} there is a positive integer $e$ such that for all $x \in R$, $\bar \nu_I(x) \in \frac{1}{e}\nn$, and
\item \label{integral-closure-I-k} if $k \geq 0$ is an integer then $\bar \nu_I(x) \geq k$ if and only if $x \in \bar{I^k}$, where $\bar{I^k}$ is the integral closure of $I^k$ in $R$.
\item \label{associated-semidegrees} Assume in addition that $I$ is a principal ideal generated by $u$ and $\bar R$ is an integral extension of $R$ which is a Krull domain. Let $\ppp_1, \ldots, \ppp_r$ be the height $1$ prime ideals of $\bar R$ containing $u$. Then for all $x \in R$, $\bar \nu_I(x) = \min\{\frac{\nu_i(x)}{e_i} : i = 1, \ldots, r\}$, where for each $i = 1, \ldots, r$, $\nu_i$ is the valuation associated with the discrete valuation ring $\bar R_{\ppp_i}$ and $e_i := \nu_i(u)$. 
\end{asparaenum}
\end{thm}

Let $\delta$ be a finitely generated subdegree on a $\kk$-algebra $A$ with a minimal presentation $\delta = \max_{1 \leq i \leq N}\delta_i$. It is not hard to see that the Noetherian-ness of $\adelta$ implies that $\delta$ is {\em integer-valued}, i.e.\ there is no non-zero $f \in A$ such that $\delta(f) = - \infty$. On the other hand, it is not obvious (e.g.\ from identity \eqref{defn-delta_i-from-f_i}) whether $\delta_i$'s can take $-\infty$ as a value or not. But the following corollary states that in fact they can {\em not}.

\begin{prop} \label{finite-prop}
Let $A$ and $\delta$ be as above. Let $B$ be any Krull domain which is also an integral extension of $\adelta$ and $\ppp_1, \ldots, \ppp_r$ be the height one primes of $B$ containing $(1)_1$. For each $j$, $1 \leq j \leq r$, define a function $\hat\delta_j$ on $A\setminus\{0\}$ by
\begin{align*}
\hat\delta_j(f) &= \delta(f) - \frac{\nu_j((f)_{\delta(f)})}{e_j} 
\end{align*}
where $\nu_j$ is the discrete valuation of the discrete valuation ring $B_{\ppp_j}$ and $e_j := \nu_j((1)_1)$. Then for each $i$, $1 \leq i \leq N$, $\delta_i \equiv \hat\delta_j$ for some $j$, $1 \leq j \leq r$. In particular, for every $i$, the semidegree $\delta_i$ is integer valued, and consequently $-\delta_i$ is a discrete valuation.
\end{prop}

\begin{proof}[Sketch of a proof]
Let $I$ be the ideal generated by $(1)_1$ in $\adelta$ and $\bar \nu_I$ be the Rees' valuation corresponding to $I$. Let $f \in A$. Since $\delta(f) \in \zz$ (from the remark following Theorem \ref{rees-thm}) and $\delta$ is a subdegree, it follows that $\bar \nu_I((f)_{\delta(f)}) = 0$. Assertion \ref{associated-semidegrees} of Rees' theorem then implies that $\min_{j=1}^r \frac{\nu_j((f)_{\delta(f)})}{e_j} = 0$, where $e_j := \nu_j((1)_1)$ for every $j$. Therefore $\delta(f) = \delta(f) - \min_{j=1}^r \frac{\nu_j((f)_{\delta(f)})}{e_j} = \max_{j=1}^r \hat\delta_j(f)$ for all $f \in A$. The assertions of Proposition \ref{finite-prop} then follow from the uniqueness of the minimal presentation of subdegrees (assertion \ref{sgf-unique-decomp} of Theorem \ref{gengfsgf-characterization}) provided we can show that each $\hat\delta_j$ is a semidegree. The latter assertion follows from case by case computations - we refer the reader to \cite[Theorem 2.2.11]{pisis} for details.
\end{proof}

Let $X$ be an affine algebraic variety and $\delta$ be a finitely generated subdegree on $A := \kk[X]$ with the minimal presentation $\delta = \max_{1 \leq i \leq N}\delta_i$. Pick $i$ such that $\delta_i \not\equiv 0$. It follows from (the proof of) Corollary \ref{sgf-components-at-infinity} that $\delta_i$ corresponds to a component $V_i$ of the hypersurface at infinity on $\xdelta$. Using Proposition \ref{finite-prop}, the following proposition establishes the relation between $\delta_i$ and the orders of vanishing along $V_i$. 

\begin{prop}\label{pole-and-degree} 
Let $X$, $A$, $\delta$, $\delta_i$ and $V_i$ be as above. For every $f \in \kk(X)\setminus\{0\}$, the order of vanishing of $f$ along $V_i$ is $-\delta_i(f)/d_i$, where $d_i := \gcd\{\delta_i(g): g \in \kk(X),\ \delta_i(g) > 0\}$. In particular, the local ring $\sheaf_{V_i, \xdelta}$ of $\xdelta$ at $V_i$ is regular. 
\end{prop}

\begin{proof}
Let $\ppp_i$ be the prime ideal of $\adelta$ corresponding to $\delta_i$. Then $\locali$ is the degree zero part of the local ring $\adelta_{\ppp_i}$. As in Remark \ref{t-remark}, let us identify $\adelta$ with $\sum F_dt^d$. Then 
\begin{align}
\begin{split}\label{local-i-expression1}
\locali &= \{\frac{ft^{kd}}{(gt^d)^k} : gt^d \not\in \ppp_i,\ d \geq \delta(g),\ k \geq 0\} \\
		&= \{\frac{f}{g^k} : (g)_{\delta(g)} \not\in \ppp_i,\ \delta(f) \leq k\delta(g),\ k \geq 0\}.
\end{split}
\end{align}
Recall from Proposition \ref{finite-prop} that $\nu_i := -\delta_i/d_i$ is a discrete valuation on $\kk(X)$. Let $R_i \subseteq \kk(X)$ be the valuation ring of $\nu_i$. It suffices to show that $\locali = R_i$. Recall that $R_i := \{g_1/g_2: g_1, g_2 \in A,\ g_2 \neq 0,\ \nu_i(g_1/g_2) \geq 0\}$. Pick $g_1, g_2 \in A$ such that $g_1/g_2 \in R_i$. Then $\nu_i(g_1/g_2) = \nu_i(g_1) - \nu_i(g_2) \geq 0$ and hence $\delta_i(g_1) \leq \delta_i(g_2)$. Pick $f_i \in A$ such that $\delta_i(f_i) > \delta_j(f_i)$ for all $j \neq i$. It follows due to \eqref{defn-delta_i-from-f_i} that there is $k \geq 1$ such that $\delta(g_lf_i^k) = \delta_i(g_lf_i^k)$ for $l = 1, 2$. Then $\delta(g_1f_i^k) = \delta_i(g_1f_i^k) = \delta_i(g_1) + \delta_i(f_i^k) \leq \delta_i(g_2) + \delta_i(f_i^k) = \delta_i(g_2f_i^k) = \delta(g_2f_i^k)$. Moreover, Lemma \ref{gensgf-characterize-lemma} implies that $(g_lf_i^k)_{\delta(g_lf_i^k)} \not\in \ppp_i$ for $l = 1, 2$. Then, according to \eqref{local-i-expression1}, $g_1f_i^k/(g_2f_i^k) = g_1/g_2 \in \locali$. Consequently $R_i \subseteq \locali$. Therefore $\locali = R_i$ due to

\begin{prolemma} \label{dvr-inclusion-maximality}
Let $R$ be a discrete valuation ring and $K$ be the quotient field of $R$. If $S$ is a proper subring of $K$ such that $R \subseteq S$, then $R = S$.
\end{prolemma}

\begin{proof}
Let $\nu$ be the discrete valuation associated to $R$ and $h \in R$ be a parameter for $\nu$, in particular $\nu(h) = 1$. Assume contrary to the claim that $R \neq S$. Let $f \in S\setminus R$. Then $f = u/h^k$ for some unit $u$ of $R$ and $k > 0$. It follows that $h^{-1} = u^{-1}fh^{k-1} \in S$. Let $g \in K\setminus \{0\}$. Then $\nu(gh^{-\nu(g)}) = 0$ and therefore $gh^{-\nu(g)} \in R \subseteq S$ and also $g = gh^{-\nu(g)} \cdot h^{\nu(g)} \in S$. Therefore $S = K$, contrary to the assumptions, which completes the proof.
\end{proof}

To summarize, we have proved that $\locali$ is precisely the valuation ring of $\nu_i$ (and therefore regular). Since valuations are completely determined by their valuation rings \cite[Section VI.8]{zsII}, it follows that $\nu_i$ is the valuation corresponding to $\locali$, which completes the proof of the proposition.
\end{proof}

\begin{rem}
If $V$ is a codimension one irreducible subvariety a variety $Y$, then $\sheaf_{V,Y}$ is a regular local ring iff $V \not\subseteq \sing Y$ iff codimension of $\sing Y \cap V$ in $Y$ is at least $2$ (where $\sing Y$ is the set of singular points of $Y$). Therefore, an equivalent formulation of Proposition \ref{pole-and-degree} is to say that in a completion $\xdelta$ of an affine variety $X$ determined by a subdegree $\delta$, the codimension of the `singular points at infinity' is at least two, or in other words, $\xdelta$ is  {\em non-singular in codimension one at infinity}.
\end{rem} 

\begin{example} \label{toric-pole-and-degree-example}
Let $\scrP$ be a convex polytope of dimension $n$ containing the origin in its interior and $\scrQ$ be a facet of $\scrP$. Recall from Example \ref{toric-example} that the semidegree corresponding to $\scrQ$ associated with the subdegree $e\delta_\scrP$ (where $e$ is an appropriate integer to ensure $e\delta_\scrP$ is integer valued) is $e\delta_\scrQ$, where $\delta_\scrQ(x^\alpha) := \frac{\langle \omega_\scrQ, \alpha \rangle}{c_\scrQ}$.
Since the greatest common divisor of the coordinates of $\omega_\scrQ$ is $1$, Proposition \ref{pole-and-degree} implies the familiar fact that the order of zero of $x^\alpha$ along the component of the hypersurface at infinity corresponding to $\scrQ$ is $- \langle \omega_\scrQ, \alpha \rangle$ \cite[Section 3.3]{fultoric}.
\end{example}

Let $X$ and $\delta$ be as in proposition \ref{pole-and-degree}. If $X$ is normal, then Proposition \ref{pole-and-degree} implies that $\xdelta$ is non-singular in codimension one, i.e.\ $\xdelta$ satisfies one of the two criteria of Serre for normality (see, e.g.\  \cite[Theorem 39]{matsulgebra}). It is not hard to see that if $X$ is normal, then $\xdelta$ is in fact normal (it follows from the fact that in this case $\kxdelta$ is integrally closed - see \cite[Proposition 2.2.7]{pisis}). If on the other hand $X$ is not normal, then clearly $\xdelta$ is not normal. But we will see below that all is not lost: $\xdelta$ is {\em relatively normal at infinity with respect to} $X$.

\begin{defn} \label{relative-normality-defn}
Let $Y$ be an algebraic variety containing $X$ as a (Zariski) dense open subset and $\psi: Z \to Y$ be a morphism of algebraic varieties. We say that $\psi$ is {\em relatively normal at infinity with respect to $X$} if for any open subset $U$ of $Y$, $\Gamma(\psi^{-1}(U), \sheaf_Z)$ is integrally closed in $\Gamma(\psi^{-1}(U \cap X), \sheaf_Z)$. In the case that $Z = Y$ and $\psi$ is the identity, we simply say that $Y$ is relatively normal at infinity with respect to $X$.
\end{defn}

\begin{prop} \label{relatively-normal-subdegree}
Let $\delta$ be a finitely generated subdegree on the ring of regular functions on an affine variety $X$. Then $\xdelta$ is relatively normal at infinity with respect to $X$.
\end{prop}

\begin{proof}
Let $A := \kk[X]$ and $D(G) := \{Q \in \xdelta | G \not \in Q\} = \spec (\adelta)_{(G)}$ be a basic affine open subset of $\xdelta$, where $G$ is a homogeneous element in \profing of {\em positive} degree $d$ and $(\profing)_{(G)}$ is the subring of $(\profing)_{G} := \profing[\frac{1}{G}]$ consisting of degree zero homogeneous elements. Then $G = (g)_d$ for some $g \in A$ with $d \geq \delta(g)$ and therefore $(\adelta)_{(G)} = \{\frac{f}{g^k}: f \in A$, $\delta(f) \leq kd\}$. In particular, $\frac{1}{g} \in (\adelta)_{(G)}$. Consequently, the regular functions on $D(G) \cap X$ are generated as a $\kk$-algebra by $A$ and $\frac{1}{g}$, i.e.\ $D(G) \cap X = \spec A_g$. We now show that $(\adelta)_{(G)}$ is integrally closed relative to $A_g$.\\

If $d > \delta(g)$, then $G = (g)_d((1)_1)^{d - \delta(g)}$, so that $\spec (\adelta)_{(G)} \subseteq \spec (\adelta)_{((1)_1)} = \spec A = X$ (Proposition \ref{prop-spec-completion}). It follows that $(\adelta)_{(G)} = A_g$, so that $(\adelta)_{(G)}$ is obviously integrally closed relative to $A_g$. So assume $d = \delta(g)$. Let $\delta_1, \ldots, \delta_N$ be the semidegrees associated to $\delta$. W.l.o.g. we may assume that there exists $m$, $1 \leq m \leq N$, such that $d = \delta_1(g) = \cdots = \delta_m(g) > \delta_i(g)$ for all $i$, $m < i \leq N$.

\begin{claim} \label{local-claim}
$(\adelta)_{(G)} = A_g \cap \bigcap_{i=1}^m \sheaf_{V_i, \xdelta}$, where for every $i$, $1 \leq i \leq m$, $V_i$ is the component of the hypersurface at infinity of $\xdelta$ corresponding to $\delta_i$.
\end{claim}

\begin{proof}
At first note that
\begin{align*} 
(\adelta)_{(G)} &= \{\frac{(f)_{kd}}{((g)_d)^k}: \delta(f) \leq kd,\ k \geq 0\} 
	    = \{\frac{f}{g^k} \in A_g: \delta(f) \leq kd,\ k \geq 0\} \\
	    &= A_g \cap \bigcap_{i=1}^N \{\frac{f}{g^k}: f \in A,\ \delta_i(f) \leq kd,\ k \geq 0\} \displaybreak[0]\\
	    &= A_g \cap \bigcap_{i=1}^m \{\frac{f}{g^k}: f \in A,\ \delta_i(f) \leq k\delta_i(g),\ k \geq 0\} 
	    \cap \bigcap_{i=m+1}^N \{\frac{f}{g^k}: f \in A,\ \delta_i(f) \leq kd,\ k \geq 0\}  \displaybreak[0] \\
	    &= A_g \cap \bigcap_{i=1}^m \sheaf_{V_i, \xdelta} \cap \bigcap_{i=m+1}^N \{\frac{f}{g^k}: f \in A,\ \delta_i(f) \leq kd,\ k \geq 0\},
\end{align*}
where the last equality is a consequence of Proposition \ref{pole-and-degree}. It follows that $(\adelta)_{(G)} \subseteq A_g \cap \bigcap_{i=1}^m \sheaf_{V_i, \xdelta}$. It remains to show the opposite inclusion. Let $f \in A$ and $\frac{f}{g^k} \in \bigcap_{i=1}^m \sheaf_{V_i, \xdelta}$. Then $\delta_i(f) \leq k\delta_i(g) = kd$ for all $i$, $1 \leq i \leq m$. Recall that $\delta_j(g) < d$ for all $j$ with $m + 1 \leq j \leq N$. It follows that if $l$ is a sufficiently large integer, then for all $j$, $m + 1 \leq j \leq N$, $\delta_j(fg^l) = \delta_j(f) + l\delta_j(g) \leq (l+k)d$. Pick an integer $l$ as in the preceding sentence. Then $\delta(fg^l) = \max_{j =1}^N \delta_j(fg^l) \leq (l+k)d$, and therefore $\frac{f}{g^k} = \frac{fg^l}{g^{k+l}} \in (\adelta)_{(G)}$. This implies that $(\adelta)_{(G)} \supseteq A_g \cap \bigcap_{i=1}^m \sheaf_{V_i, \xdelta}$ and completes the proof of the claim.
\end{proof}

Recall (from Proposition \ref{pole-and-degree}) that each $\sheaf_{V_i, \xdelta}$ is a discrete valuation ring, therefore in particular is integrally closed. It then follows that $(\adelta)_{(G)}$ is integrally closed relative to $A_g$, as required to prove the proposition.
\end{proof}

\begin{rem}
In the same way that normality implies non-singularity in codimension one, it can be shown that for a completion $\bar X$ of $X$, if $\bar X$ is normal at infinity with respect to $X$, then $\bar X$ is also non-singular at infinity in codimension one (with respect to $X$). Therefore Proposition \ref{relatively-normal-subdegree} strengthens Proposition \ref{pole-and-degree}.
\end{rem}

In general, completions corresponding to degree-like functions may have arbitrarily bad singularities at infinity. Below we introduce the notion of {\em normalization at infinity}, which makes the hypersurface at infinity of a given completion a bit `less singular', in the same sense that normalization of a singular variety is less singular than itself.  

\begin{defn}
Let $Y$ be a variety containing $X$ as a dense open subset. The {\em normalization of $Y$ at infinity with respect to $X$} is another variety $\tilde Y$ containing $X$ as a dense open subset and a finite morphism $\phi: \tilde Y \to Y$ such that $\tilde Y$ is normal at infinity with respect to $X$ and  $\phi|_X$ is the identity map.
\end{defn}

\begin{prop} \label{normalization-at-infinity}
Let $X$ be an arbitrary irreducible variety and $Y$ be a variety containing $X$ as a Zariski dense open subset. Then the normalization $\tilde Y$ of $Y$ at infinity (with respect to $X$) exists and is unique. Moreover, $\tilde Y$ has the following universal properties: 
\begin{compactenum}
\item \label{normal-universal-1}If $\psi: Z \to Y$ is a dominant morphism of algebraic varieties such that $Z$ is normal at infinity with respect to $X$ via $\psi$, then there exists a unique morphism $\theta: Z \to \tilde Y$ such that the following diagram commutes.
$$\xymatrix{
 			& Z \ar[1,-1]_\theta \ar[1,1]^\psi 	& \\
\tilde Y \ar[0,2]^\phi	& 					& Y}$$
\item \label{normal-universal-2} If $Z$ is another variety containing $X$ as an open subset and $\psi: Z \to Y$ is a finite morphism such that $\psi|_X$ is the identity map, then there exists a unique morphism $\theta: \tilde Y \to Z$ such that the following diagram commutes.
$$\xymatrix{
 					  & Z \ar[1,1]^\psi 	& \\
\tilde Y \ar[-1,1]^\theta \ar[0,2]^\phi &			& Y}$$
\end{compactenum}
\end{prop}

\begin{proof}
The universal properties of the normalization at infinity are completely analogous to those of normalization and the standard proofs of the latter (e.g.\ the proof in \cite{shaf1}) apply almost word by word to prove the former. We therefore skip the details. 
\end{proof}

Let $X$ be a Zariski open subset of a variety $Y$. If $Y$ is affine (resp.\ projective), then it can be shown that the normalization $\tilde Y$ of $Y$ at infinity with respect to $X$ is also affine (resp.\ projective). If $X$ is affine and $Y$ is a completion of $X$ determined by a degree-like function $\delta$ on $\kk[X]$, then similarly it turns out that $\tilde Y$ is also the completion of $X$ determined by a degree-like function $\tilde \delta$. In fact, as the next theorem shows, $\tilde \delta$ can be taken to be a {\em subdegree}. 

\begin{thm}[Main Existence Theorem]\label{thm-noeth-integral-subdegree}
Let $X$ be an affine variety and $\delta$ be a finitely generated degree-like function on $A := \kk[X]$. Then 
\begin{enumerate}
\item \label{normalization-defn} There is a positive integer $e$ and a subdegree $\tilde\delta$ on $A$ such that for all $h \in  A$, 
$$\tilde \delta(h) := e\lim_{m \to \infty} \frac{\delta(h^m)}{m}.$$
\item \label{normalization-finitely-gen} There is a natural inclusion $A^{e\delta} \subseteq A^{\tilde \delta}$ of graded rings such that $A^{\tilde \delta}$ is integral over $A^{e\delta}$. In particular, $\tilde \delta$ is finitely generated. If $\delta$ is non-negative (resp.\ projective), then $\tilde \delta$ is also non-negative (resp.\ projective).
\item \label{normalizing-subdegree-normalizes} The variety $\xtildedelta$ is the normalization at infinity of $\xdelta$ with respect to $X$.
\end{enumerate}
\end{thm}

\begin{proof}
Let $I$ be the ideal generated by $(1)_1$ in $\adelta$. Fix $h \in A$ and $m \in \nn$. Since $\delta(h^m) \leq m\delta(h)$, it follows that $k:=m\delta(h) - \delta(h^m)$ is the largest integer such that $(h^m)_{m\delta(h)} \in I^{k}$, and consequently $\nu_I(((h)_{\delta(h)})^m) = k = m\delta(h) - \delta(h^m)$ (where $\nu_I$ is defined as in \eqref{nu-nu-bar}). Therefore $\delta(h^m)/m = \delta(h) - \nu_I(((h)_{\delta(h)})^m)/m$. Assertion \ref{normalization} of Theorem \ref{rees-thm} then implies that $\bar\delta(h) := \lim_{m \to \infty} \delta(h^m)/m = \delta(h) - \bar \nu_I((h)_{\delta(h)})$ is well defined and there exists a positive integer $e$ (independent of $h$) such that $\bar\delta(h) \in \frac{1}{e}\zz$. Taking $\tilde \delta := e\bar \delta$ proves the displayed formula of assertion \ref{normalization-defn}. Moreover, note that $\tilde \delta(h^m) = m\tilde \delta(h)$ for all $h$ and $m$. Therefore, if we show that $A^{\tilde \delta}$ is finitely generated, then it follows via Corollary \ref{sgf-reformulation} that $\tilde\delta$ is a subdegree. Hence it suffices to prove assertion \ref{normalization-finitely-gen} to complete the proof of assertion \ref{normalization-defn}. \\

We now prove assertion \ref{normalization-finitely-gen}. Let $\scrF := \{F_d\}_{d \in \zz}$ be the filtration on $A$ corresponding to $\delta$ and identify $\adelta$ with $\bigoplus_{i \geq 0} F_it^i$. Let $\bar F_{\frac{m}{e}} := \{h \in A: \bar\delta(h) \leq \frac{m}{e}\}$ for all $m \geq 0$, and define $\profinggg{\bar\delta} := \bigoplus_{m \geq 0} \bar F_{\frac{m}{e}}t^{\frac{m}{e}}$. Since $\bar \delta \leq \delta$, it follows that $F_k \subseteq \bar F_k$ for each $k \in \zz$. Therefore $\adelta \subseteq \profinggg{\bar\delta}$. 

\begin{claim} \label{integrality-in-main}
$\profinggg{\bar\delta}$ is integral over $\adelta$.
\end{claim}

\begin{proof}
It suffices to show that $ht^{\bar\delta(h)}$ is integral over $\adelta$ for each $h \in  A$ such that $\bar \delta(h) \geq 0$. Let $h$ be as in the preceding sentence. Then $ht^{\bar\delta(h)}$ is integral over $\adelta$ if and only if $\bar H := (ht^{\bar\delta(h)})^e $ is integral over $\adelta$. Note that $e\bar\delta(h) = \bar\delta(h^e)$ by construction of $\bar\delta$, and therefore $\bar H = h^et^{\bar\delta(h^e)}$. Let $H := h^et^{\delta(h^e)} \in \adelta$ and $k := \bar \nu_I(H)$, where $I$ is the ideal generated by $(1)_1$ in $\adelta$. Since $\nu_I(H^m) = m\delta(H) - \delta(H^m)$ and $\delta(H) = \delta(h^e)$, it follows that $k =  \delta(h^e) - \bar\delta(h^e) = \delta(h^e) - e\bar\delta(h)$. Then $k$ is an integer. Hence according to assertion \ref{integral-closure-I-k} of Theorem \ref{rees-thm}, $H$ is in the integral closure of $I^{k}$ in $\adelta$, i.e.\ $H$ satisfies an equation of the form $H^l + G_1H^{l-1}+ \cdots + G_l = 0$, where $G_i \in I^{ik}$ for each $i$. Comparing the coefficients at $t^{l\delta(h^e)}$ in the above equation, we may assume w.l.o.g.\ that each $G_i$ is of the form $g_it^{i\delta(h^e)}$ for some $g_i \in A$. Since $G_i \in I^{ik}$, it follows that $i\delta(h^e) \geq \delta(g_i) + ik$, implying that $g_it^{i(\delta(h^e) - k)}$ is an element of $\adelta$. Moreover, regarding $H$ as an element of the ring $\profinggg{\bar\delta}$ (via the embedding $\adelta \into \profinggg{\bar\delta}$) yields $H = h^et^{\delta(h^e)} = h^et^{\bar\delta(h^e)+k} = h^et^{\bar\delta(h^e)}t^k = t^k\bar H$. Substituting these values of $H$ and $G_i$ into the equation of integral dependence for $H$ and then canceling a factor of $t^{lk}$ we conclude that $(\bar H)^l + \sum_{i=1}^l g_it^{i(\delta(h^e) - k)}(\bar H)^{l-i} = 0$. But then $\bar H$ is integral over $\adelta$, which completes the proof of the claim.
\end{proof}

Let $\tilde \scrF := \{\tilde F_d\}_{d \geq 0}$ be the filtration corresponding to $\tilde \delta = e \bar \delta$. Observe that $\tilde F_d = \{f : e\bar\delta(h) \leq d\} = \bar F_{\frac{d}{e}}$, so that the homomorphism $\chi: \profinggg{\bar\delta} := \bigoplus_{d \geq 0} \bar F_{\frac{d}{e}}t^{\frac{d}{e}} \longrightarrow \bigoplus_{d \geq 0}\tilde F_ds^d \cong \profinggg{\tilde\delta}$ that sends $t \mapsto s^e$ and is the identity map on the coefficients (i.e.\ on $\bar F_{\frac{d}{e}}$ for $d \geq 0$) is in fact an isomorphism of $\kk$-algebras. Therefore, it follows due to Claim \ref{integrality-in-main} that $A^{\tilde \delta} = \chi(A^{\bar \delta})$ is integral over $\chi(\adelta)$. On the other hand, since $\tilde \delta \leq e\delta$, there is a natural inclusion $\profinggg{e\delta} \subseteq \profinggg{\tilde\delta}$ of graded rings and $\chi(\adelta) \subseteq \profinggg{e\delta}$. Therefore $\profinggg{\tilde\delta}$ is integral over $\profinggg{e\delta}$. Since $\adelta$ (and therefore also $A^{e\delta}$) is a finitely generated $\kk$-algebra, it follows that $\profinggg{\tilde\delta}$ is a finitely generated $\kk$-algebra.\\

If $\delta$ is non-negative (resp. complete), then by construction $\tilde \delta$ is also non-negative (resp. complete). This completes the proof of assertion \ref{normalization-finitely-gen} and therefore also assertion \ref{normalization-defn} of the theorem. Moreover, Proposition \ref{relatively-normal-subdegree} implies that $\xtildedelta$ is normal at infinity with respect to $X$. Since the natural morphism $\xtildedelta \to \xdelta$ induced by the integral inclusion $\profinggg{e\delta} \subseteq \profinggg{\tilde\delta}$ is finite, it follows that $\xtildedelta$ is the normalization of $\xdelta$ at infinity with respect to $X$. This completes the proof of the theorem.
\end{proof}

\begin{reminition}\label{normalizing-subdegree}
Let $X, A, \delta$ and $\tilde \delta$ be as in Theorem \ref{thm-noeth-integral-subdegree}. If in addition $X$ is normal, then $X^{\tilde \delta}$ is also normal (cf.\ the remark following Example \ref{toric-pole-and-degree-example}). The universal property of the normalization then implies that $\xtildedelta$ is in fact the normalization of $\xdelta$. Motivated by this, we will refer to $\tilde \delta$ as a {\em normalization} of $\delta$. (Note that the normalization of degree-like functions is unique only up to a constant factor, depending on the choice of $e$ of assertion 1 of Theorem \ref{thm-noeth-integral-subdegree}).
\end{reminition}

\begin{example}
Let $\scrA$ be a finite subset of $\zz^n$ such that $\zz^n = \zz\scrA$ and the convex hull $\scrP$ of $\scrA$ in $\rr^n$ contains the origin in its interior. Let $X_\scrA$ be the closure of the image of the map $\phi_\scrA: \nktorus \into \pp^{|\scrA| -1}(\kk)$ whose components are the monomials $x^\alpha$ with $\alpha \in \scrA$. $X_\scrA$ is the (possibly non-normal) toric variety corresponding to $\scrA$ (see e.g.\ \cite[Section 5.1]{gkz}). Let $X := \nktorus$ and $\eta$ be the degree-like function on $\kk[X]$ corresponding to the completion $X \into X_\scrA$ and $e\delta_\scrP$ be a subdegree associated to $\scrP$ as in Example \ref{toric-example}. Then it is not hard to see that $\delta_\scrP$ is a normalization of $\eta$. This implies (following Remark-Definition \ref{normalizing-subdegree}) that $X_\scrP$ is the normalization of $X_\scrA$ (cf.\ \cite[Proposition 2.8(a)]{gkz}).
\end{example}

\begin{rem} \label{normalizing-remark}
Let $\delta$ be a finitely generated degree-like function on $\kk[X]$ and $\tilde \delta$ be a normalization of $\delta$. Theorem \ref{thm-noeth-integral-subdegree} states that applying $\tilde \delta$ to $\kk[X]$ produces the normalization at infinity of $\xdelta$ with respect to $X$. In light of this statement, it is natural to ask if applying $\tilde \delta$ to the coordinate ring of the normalization of $X$ would produce the normalization of $\xdelta$. This indeed turns out to be true, i.e.\ if we denote normalization by the symbol `tilde', then we have: 
\begin{align*}
\bar{\tilde X}^{\tilde \delta} \cong \widetilde{\xdelta}. 
\end{align*}
\end{rem}

We end this section with an application of our results. Namely, we prove the following `finiteness-property' of divisorial valuations (defined in the paragraph following Theorem \ref{weak-structure}) proved in \cite{fernex-ein-ishii} for $\kk = \cc$.

\begin{thm} \label{finiteness-prop}
Let $X$ be an irreducible affine variety and $\nu$ be a divisorial valuation over $X$. Then there exist elements $f_1, \ldots, f_r \in \kk[X]\setminus\{0\}$ such that for every $f \in \kk[X] \setminus\{0\}$,
\begin{align}
\nu(f) &= \min\{\nu'(f) : \nu'\ \text{is a divisorial valuation over $X$,}\ \nu'(f_i) = \nu(f_i)\  \text{for}\ 1 \leq i \leq r \} \tag{$*$} \\
	&= \min\{\nu'(f) : \nu'\ \text{is a valuation on $\kk(X)$ such that the value group of $\nu$ contains} \notag \\
	& \mbox{\phantom{=$\min\{\nu'(f):$}}\ \ \text{the integers and}\ \nu'(f_i) = \nu(f_i)\ \text{for}\ 1 \leq i \leq r \} \tag{$*'$} 
\end{align}
\end{thm}

\begin{proof}
There exists a normal variety $Y$ equipped with a birational morphism $\pi: Y \to X$ and a prime divisor $E$ on $Y$ such that $\nu = q\ord_E$ for some positive integer $q$. Let $U$ be an open affine subset of $Y$ such that $E \cap U \neq \emptyset$ and $\bar U$ be an arbitrary projective completion of $U$. Pick an effective ample divisor $D$ on $\bar U$ such that $E \cap U \subseteq \supp(D)$ and set $Z := \bar U \setminus \supp(D)$. Then $\bar U \cong \bar Z^\delta$ for some degree-like function $\delta$ on $\kk[Z]$ (Proposition \ref{prop-spec-completion}). Let $\tilde \delta$ be a normalization of $\delta$ with minimal presentation $\tilde \delta = \max_{i=1}^N \delta_i$. Proposition \ref{pole-and-degree} implies that there exists $i$ such that $\delta_i = -\frac{d}{q}\nu$ for some positive integer $d$. W.l.o.g.\ assume $i = 1$. Pick $g_1 \in \kk[Z]\setminus\{0\}$ such that $\tilde \delta(g_1) = \delta_1(g_1) > \delta_i(g_1)$ for all $i > 1$. Finally, pick $g_2, \ldots, g_s \in \kk[Z]\setminus\{0\}$ such that $(g_1)_{\tilde \delta(g_1)}, \ldots, (g_s)_{\tilde \delta(g_s)}$ generate $\kk[Z]^{\tilde \delta}$ as a $\kk$-algebra. We claim that for all $g \in \kk[Z]\setminus\{0\}$,
\begin{align}
\delta_1(g) &= \max\{\delta'(g) : \delta'\ \text{is a semidegree on $\kk[Z]$,}\ \delta'(g_i) = \delta_1(g_i)\ \text{for}\ 1 \leq i \leq s \} \tag{$\tilde{*}$}. \label{finiteness-tilde}
\end{align}
Indeed, pick a semidegree $\delta'$ on $\kk[Z]$ such that $\delta'(g_i) = \delta_1(g_i)$ for all $i$, $1 \leq i \leq s$. Since $\delta_1 \leq \tilde \delta$, it follows from the choice of $g_j$'s that $\delta' \leq \tilde \delta$. Now pick any $g \in \kk[Z]\setminus\{0\}$. Then for all sufficiently large $k$, we have that $\delta_1(gg_1^k) = \tilde \delta(gg_1^k) \geq \delta'(gg_1^k)$, and consequently, $\delta_1(g) \geq \delta'(g)$. This proves identity \eqref{finiteness-tilde}. \\

Let $g_0 \in \kk[Z] \setminus\{0\}$ such that $Z \setminus V(g_0) \subseteq X$ and let $f_1, \ldots, f_r \in \kk[X]\setminus\{0\}$ such that for each $i$, $0 \leq i \leq s$, $g_i = f_j/f_k$ for some $j, k$. We show that $f_1, \ldots, f_r$ satisfies the claim of the theorem. Indeed, pick any $\nu'$ as in the right hand side of the identity \eqref{finiteness'} and an arbitrary $f \in \kk[X] \setminus \{0\}$. Then $f = g/g_0^k$ for some $k \geq 0$ and $g \in \kk[Z]$. Since $\nu'(g_j) = \nu(g_j) = -\frac{q}{d}\delta_1(g_j)$ for all $j$, $0 \leq j \leq s$, applying identity \eqref{finiteness-tilde} to $\delta' := -\frac{d}{q}\nu'$ yields that $\nu(f) \leq \nu'(f)$, as required.
\end{proof}

\section{Divisor at infinity} \label{sec-sub-divisor}
Associated to every degree-like function $\delta$ there is a canonical ample $\qq$-Cartier divisor supported at infinity on the corresponding completion. In this section we study this divisor for the case that $\delta$ is a subdegree and establish a formula for its pull-back under a dominant morphism. As an application we compute the matrix of intersection numbers of the curves at infinity on a class of completions of certain affine surfaces. We also compute the nef cone of the latter completions and as a consequence give a positive answer to Question \ref{finite-question} in the case of surfaces under the additional hypothesis that each $\delta_j$ is strictly positive on every non-constant regular function on $X$ (Corollary \ref{finite-corollary}).

\subsection{Pull-back formula}

\newcommand\divinfinity[1]{%
	D^{#1}_{\infty}%
}
\newcommand\divdeltainfinity{%
	\divinfinity{\delta}%
}
\newcommand\divetainfinity[1]{%
	\divinfinity{\eta}%
}

\newcommand\xinfinity[1]{X^{#1}_\infty}
\newcommand\xdeltainfinity{\xinfinity{\delta}}

\begin{reminition}
Let $X$ be an affine variety and $\delta$ be a finitely generated degree-like function on $\kk[X]$. The {\em divisor at infinity} $\divdeltainfinity$ on $\xdelta$ is the $\qq$-Cartier divisor defined by $(1)_1$. $\divdeltainfinity$ is Cartier iff there is a subset $S$ of $\kk[X]^\delta$ consisting of homogeneous elements of degree one such that the radical of the ideal generated by $S$ is the {\em irrelevant ideal} of $\kk[X]^\delta$. It follows from Proposition \ref{prop-spec-completion} that $\divdeltainfinity$ is ample and supported at {\em infinity}, i.e. $\supp \divdeltainfinity = \xdelta \setminus X =: \xdeltainfinity$.
\end{reminition}

\begin{notation*}
Let $Z$ be an algebraic variety. For a Cartier divisor $D$ on $Z$ (resp.\ an irreducible codimension one subvariety $V$ of $Z$), we denote the corresponding Weil divisor by $[D]$ (resp.\ $[V]$). Moreover, given a $\qq$-Weil or $\qq$-Cartier divisor $D$ on $Z$, we write $D \geq 0$ if $D'$ is {\em effective}. For convenience of the reader we recall the definition of `effectiveness' of $\qq$-Weil and $\qq$-Cartier divisors: if $D$ is a $\qq$-Weil divisor, then $D$ is effective iff for each irreducible codimension one subvariety $V$ of $Z$, the coefficient of $[V]$ in $D$ is zero. On the other hand, if $D$ is Cartier, then $D$ is said to be effective iff all the local equations of $D$ are regular functions. Finally, a $\qq$-Cartier divisor is effective iff it has a multiple which is an effective Cartier divisor.
\end{notation*}

\begin{lemma} \label{divisor-at-infinity}
Let $X$ be an irreducible affine variety, $\delta$ be a finitely generated subdegree on $\kk[X]$, and $\delta_1, \ldots, \delta_N$ be the non-zero semidegrees associated to $\delta$. Then the $\qq$-Weil divisor associated to $\divdeltainfinity$ is
\begin{align*}
[\divdeltainfinity] = \sum_{j=1}^N \frac{1}{d_j}[V_j]
\end{align*}
where for each $j$, $1 \leq j \leq N$, $d_j := \gcd\{\delta_j(g): g \in \kk(X),\ \delta_j(g) > 0\}$ and $V_j$ is the irreducible component of $\xdeltainfinity$ corresponding to $\delta_j$.
\end{lemma}

\begin{proof}
Fix $j$, $1 \leq j \leq N$ and $g \in \kk[X]$ such that $d := \delta(g) = \delta_j(g)$. Let $U := \xdelta \setminus V((g)_d)$. Assertion \ref{part1-gensgf-characterize-lemma} of Lemma \ref{gensgf-characterize-lemma} implies that $U \cap V_j \neq \emptyset$. Recall that for each $f \in \kk(X) \setminus\{0\}$, the order of vanishing of $f$ at $V_j$ is precisely $-\delta_j(f)/d_j$ (Proposition \ref{pole-and-degree}). Since $1/g$ is a local equation of $d\divdeltainfinity$ on $U$, it follows that the coefficient of $[V_j]$ in the expression for $[d\divdeltainfinity]$ is $-\frac{1}{d_j}\delta_j(1/g) = d/d_j$, as required to prove the lemma.
\end{proof}

\begin{example} 
Let $\scrP$ be an $n$-dimensional convex rational polytope as in Example \ref{toric-example} and $X_\scrP$ be the corresponding toric variety. We use the notation of Example \ref{toric-example}. With $X := (\kk^*)^n$ and $\delta := e\delta_\scrP$, it follows that $[\divdeltainfinity] = \sum_\scrQ \frac{c_\scrQ}{e} [X_\scrQ]$, where the sum is over all facets $\scrQ$ of $\scrP$ and $X_\scrQ$ are the components of $X_\scrP \setminus X$ corresponding to $\delta_\scrQ$ (see Example \ref{toric-number-of-components-example}). In particular, $\divdeltainfinity  = \frac{1}{e}D_\scrP$, where $D_\scrP$ is the divisor on $X_\scrP$ corresponding to $\scrP$ \cite[Section 3.4]{fultoric}.
\end{example}

Next we establish the pull-back formula for the divisor at infinity. To state it we use from \cite[Section 1.2]{fultersection} the notion of {\em order of vanishing} $\ord_W(f)$ of a non-zero rational function $f$ on an arbitrary irreducible variety $Z$ along a codimension one irreducible subvariety $W$. 

\begin{prop} \label{pulling-back-infinity}
Let $X$ be an irreducible affine variety, $\delta$ be a finitely generated non-negative (non-zero) subdegree on $\kk[X]$ and $\phi: Z \to \xdelta$ be a dominant morphism from an irreducible variety $Z$. Then
\begin{align}
[\phi^*(\divdeltainfinity)] &= \sum_W l^\phi_\infty(\delta, \pole_W) [W], \label{pull-back-formula}
\end{align}
where the sum is over codimension one irreducible subvarieties $W$ of $Z$, and for each such $W$, the function $\pole_W$ is the negative of $\ord_W$, and  
$$l^\phi_\infty(\delta, \pole_W)  := \max \left \{\frac{\pole_W(\phi^*(f))}{\delta(f)}: f \in \kk[X],\ \delta(f) > 0\right\}.$$
\end{prop}   

\begin{rem}
Recall (from Remark-Definition \ref{linking-reminition}) that $l^\phi_\infty(\delta, \pole_W)$ is the {\em linking number at infinity} ({\em relative to $\phi$}) of $\delta$ and $\pole_W$.
\end{rem}

\begin{proof}
Let $\delta_1, \ldots, \delta_N$ be the non-zero semidegrees associated to $\delta$. For each $j$, $1 \leq j \leq N$, let $V_j$ be the component of the hypersurface at infinity on $\xdelta$ corresponding to $\delta_j$ and $d_j := \gcd\{\delta_j(g): g \in \kk(X),\ \delta_j(g) > 0\}$. Pick $f \in \kk[X]$ such that $\delta(f) > 0$ and let $\Div_{\xdelta}(f)$ be the principal Cartier divisor of $f$ on $\xdelta$. Then Proposition \ref{pole-and-degree} implies that
\begin{align}
	& [\Div_{\xdelta}(f)] +  \sum_{j=1}^N \frac{\delta_j(f)}{d_j}[V_j] \geq 0 \notag \\
\im	& [\Div_{\xdelta}(f)] + \delta(f) \sum_{j=1}^N \frac{1}{d_j}[V_j] \geq 0\ \text{(since $\delta(f) \geq \delta_j(f)$ for all $j$, $1 \leq j \leq N$)}\notag \\
\im	& \Div_{\xdelta}(f) + \delta(f)\divdeltainfinity \geq 0\ \text{(due to Lemma \ref{divisor-at-infinity} and Claim \ref{local-claim}).} \label{effectively-cartier}
\end{align}
Let $W$ be an irreducible codimension one subvariety of $Z$ and $c_W$ be the coefficient of $[W]$ in $[\phi^*(\divdeltainfinity)]$. Since the pull-back (under a dominant map) of an effective $\qq$-Cartier divisor is also effective, inequality \eqref{effectively-cartier} implies that $\ord_W(\phi^*(f)) + \delta(f)c_W \geq 0$, or in other words, $c_W \geq \pole_W(\phi^*(f))/\delta(f)$. It follows that $c_W \geq l^\phi_\infty(\delta, \pole_W)$.\\

To complete the proof of the proposition, it suffices to show that there exists $f \in \kk[X]$ such that $\delta(f) > 0$ and $c_W = \pole_W(\phi^*(f))/\delta(f)$. We divide the proof into two cases: 

\paragraph{Case 1: $\phi(W) \subseteq X$.} In this case $c_W = 0$, since $\phi(W) \cap \supp \divdeltainfinity = \emptyset$. On the other hand, for all $g \in \kk[X]$, $\phi^*(g)$ is regular on a neighbourhood of $W$, so that $\pole_W(\phi^*(g)) \leq 0$. Pick any $g \in \kk[X]$ such that $\delta(g) > 0$. Then setting $f := g - \alpha$ for generic $\alpha \in \kk$ yields that $\pole_W(\phi^*(f))/\delta(f) = 0$, as required. 

\paragraph{Case 2: $\phi(W) \subsetneq X$.} In this case there exists $w \in W$ and $j$, $1 \leq j \leq N$, such that $\phi(w) \in V_j$. Since $\divdeltainfinity$ is ample, there exists $d > 0$ such that $\sheaf_{\xdelta}(d\divdeltainfinity)$ is globally generated. Note that the global sections of $\sheaf_{\xdelta}(d\divdeltainfinity)$ are precisely $\{f \in \kk[X]: \delta(f) \leq d\}$. In particular, there exists $f \in \kk[X]$ such that $\delta(f) \leq d$ and $fh$ is invertible near $\phi(w)$, where $h$ is a local equation for $d\divdeltainfinity$ near $\phi(w)$. Then $\delta_j(f) = -d_j\ord_{V_j}(f) = d_j\ord_{V_j}(h) = d$. Since $\delta(f) \geq \delta_j(f)$, it follows that $\delta(f) = d$; in particular, $\delta(f) > 0$. Moreover, since $\phi^*(fh)$ is invertible near $w$, it follows that $0 = \ord_W(\phi^*(fh)) = \ord_W(\phi^*(f)) + \ord_W(\phi^*(h)) =  -\pole_W(\phi^*(f)) + dc_W$. Taken together, the preceding two sentences imply that $c_W = \pole_W(\phi^*(f))/\delta(f)$, as required.
\end{proof}

\subsection{Applications of the pull-back formula}
Throughout this subsection we set $X$ to be an affine surface such that the only invertible regular functions on $X$ are (non-zero) constants (the most typical example being $X = \affine{2}{\kk}$). Let  $\bar X_1, \ldots, \bar X_k$ be non-isomorphic normal projective completions of $X$ such that the complement $C_j$ of $X$ in each $\bar X_j$ is an irreducible curve and let $\bar X$ be the normalization of the closure in $\bar X_1 \times \cdots \times \bar X_k$ of the image of $X$ under the {\em diagonal} mapping. As the first application of Proposition \ref{pulling-back-infinity}, we compute the matrix of intersection numbers of the curves at infinity on $\bar X$.\\

The complement of $X$ in $\bar X$ has precisely $k$ irreducible components $\tilde C_1, \ldots, \tilde C_k$, where $\tilde C_j$ is the unique curve in $\bar X \setminus X$ which maps {\em onto} $C_j$ under the natural projection $\pi_j : \bar X \to \bar X_j$. The intersection numbers $(\tilde C_i, \tilde C_j)$ are well defined (according to Mumford's intersection theory for normal complete surfaces \cite{mumford-normal}). Let $\delta_j : \kk[X] \to \nn$ be the order of pole along $C_j$. Then each $\delta_j$ is a projective semidegree. Define three $k \times k$ matrices $\scrL, \scrI, \scrD$ as follows:
\begin{align*}
\scrL &:= \text{matrix of linking numbers at infinity of $\delta_j$'s with $(i,j)$-th entry being}\ l_\infty(\delta_i, \delta_j),\\
\scrD &:= \text{the diagonal matrix with $i$-th diagonal entry being}\ (C_i, C_i),\\
\scrI &:= \text{matrix of intersection numbers $\tilde C_j$'s with $(i,j)$-th entry being}\ (\tilde C_i, \tilde C_j).
\end{align*}

\begin{lemma} \label{linking-intersections}
$\scrL \scrI = \scrD$. 
\end{lemma}

\begin{proof}
Pick $i, j$, $1 \leq i, j \leq k$. If $i \neq j$, then $\pi_i(\tilde C_j)$ is a point and therefore $(\pi_i^*(C_i), \tilde C_j) = 0$. Since intersection numbers are preserved by pull-backs under birational morphisms, it follows that $(C_i, C_i) = (\pi_i^*(C_i), \pi_i^*(C_i)) = (\pi_i^*(C_i), \tilde C_i)$. The lemma now follows from Proposition \ref{pulling-back-infinity} which implies that $\pi_i^*(C_i) = \sum_{l=1}^k l_{ik}\tilde C_k$. 
\end{proof}

As an application of Lemma \ref{linking-intersections} we compute the intersection theory of a class of toric completions of $\affine{2}{\kk}$. Let $\scrP$ be a convex rational polygon in $\rr^2$ with\\
\begin{multicols}{2}
\begin{enumerate}
\item one vertex at the origin,
\item two edges along the axes, and 
\item the other edges being line segments with {\em negative} rational slopes.  
\end{enumerate}
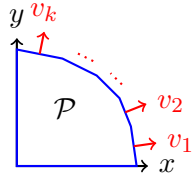
\columnbreak

\begin{tikzpicture}[scale=0.3]
	\pgfmathsetmacro\ypt{2/5 + 0.75 + 1 + 1.25 + 1.75};
	\pgfmathsetmacro\xpt{2 + 1.5 + 1 + 0.5 + 0.25};
	\draw [<->,thick] (0,\ypt + 0.5) node (yaxis) [above] {$y$}
       	 |-  (\xpt+0.5,0) node (xaxis) [right] {$x$};
	\draw [blue, thick] (0,\ypt) coordinate (p1) -- ++(2,-2/5) coordinate (p2) -- ++(1.5,-0.75) coordinate (p3) 
	-- ++(1,-1) coordinate (p4) -- ++(1*0.5,-2.5*0.5) coordinate (p5) -- ++(0.25,-7*.25) coordinate (p6)
	-- (0,0) coordinate (o) -- cycle;
	\coordinate (m1) at ($(p6)!0.5!(p5)$);
	\draw [red, ->, thick] (m1) -- ($(m1)!1cm!-90:(p5)$) coordinate (v1) node (v1label) [sloped,right] {$v_1$};
	\coordinate (m2) at ($(p5)!0.5!(p4)$);
	\draw [red, ->, thick] (m2) -- ($(m2)!1cm!-90:(p4)$) coordinate (v2) node (v2label) [sloped,right] {$v_2$};
	
	\coordinate (m31) at ($(p3)!0.25!(p4)$);
	\coordinate (m32) at ($(p3)!0.75!(p4)$);
	\draw [red, dotted, thick] ($(m32)!0.5cm!90:(p4)$) -- ($(m31)!0.5cm!90:(p4)$);

	\coordinate (m41) at ($(p3)!0.25!(p2)$);
	\coordinate (m42) at ($(p3)!0.75!(p2)$);
	\draw [red, dotted, thick] ($(m41)!0.5cm!-90:(p2)$) -- ($(m42)!0.5cm!-90:(p2)$);

	\coordinate (m5) at ($(p2)!0.5!(p1)$);
	\draw [red, ->, thick] (m5) -- ($(m5)!1cm!-90:(p1)$) coordinate (vk) node (vklabel) [sloped,above] {$v_k$};

	\draw (\xpt/2 - 0.5, \ypt/2) node (P) {$\scrP$};
\end{tikzpicture}
\end{multicols}

List the non-axis edges of $\scrP$ counter-clockwise by $e_1, \ldots, e_k$ and for each $j$, $1 \leq j \leq k$, let $v_j := (v_{j1}, v_{j2})$ be the smallest integral vector on the {\em outward pointing} normal to $e_j$ and $\tilde C_j$ be the torus invariant curve associated to $e_j$ on the toric surface $X_\scrP$ corresponding to $\scrP$. Let $\scrI$ be the $k \times k$ matrix of intersection numbers $(\tilde C_i, \tilde C_j)$ for $1 \leq i, j \leq k$. Each $v_j$ corresponds to a toric surface $\bar X_j$ corresponding to a triangle $\scrP_j$ which has two edges along the positive axes and the other edge is parallel to $e_j$. Then an application of Lemma \ref{linking-intersections} to $\bar X_1, \ldots, \bar X_k$ and $\bar X := X_\scrP$ yields:

\begin{cor} \label{toric-intersections}
$$\scrI = \left( \begin{array}{ccccc}
		1 & \frac{v_{22}}{v_{12}} & \frac{v_{32}}{v_{12}} & \cdots & \frac{v_{k2}}{v_{12}} \\
		\frac{v_{11}}{v_{21}} & 1 & \frac{v_{32}}{v_{22}} & \cdots & \frac{v_{k2}}{v_{22}} \\
		\vdots 	& \vdots & \vdots & \vdots & \vdots \\
		\frac{v_{11}}{v_{k1}} & \frac{v_{21}}{v_{k1}} & \frac{v_{31}}{v_{k1}} & \cdots & 1
		\end{array}
	\right)^{-1} 
	\left(  \begin{array}{cccc}
		\frac{1}{v_{11}v_{12}} & 0 & \cdots & 0 \\
		0 & \frac{1}{v_{21}v_{22}} & \cdots & 0 \\
		\vdots & \vdots & \ddots & \vdots\\
		0 & 0 & \cdots & \frac{1}{v_{k1}v_{k2}}
		\end{array}
	\right) \qed
$$		
\end{cor}

Let $X$ and $\bar X$ be as in the first paragraph of this subsection. Lemma \ref{linking-intersections} implies that the {\em cone $\NE(\bar X)$ of curves} on $\bar X$ is generated by (the equivalence classes) of $\tilde C_1, \ldots, \tilde C_k$ and consequently is simplicial of dimension $k$. The {\em nef cone} $\nef(\bar X)$ of $\bar X$ is dual to $\NE(\bar X)$ under the intersection product. We next give another description of $\nef(\bar X)$.

\begin{lemma} \label{nef-lemma}
$\nef(\bar X) = \{0\} \cup \{\sum_{i=1}^k m_i [\tilde C_i]: m_i > 0$ and there exist non-constant regular functions $f_i$ on $X$ such that $\delta_i(f_i)/m_i \geq \delta_j(f_i)/m_j$ for all $i,j,\ 1 \leq i,j \leq k\}$. 
\end{lemma}

\begin{proof}
Let $\scrC$ be the set of divisors on the right hand side of the equality to be proved. Since $\delta_i$'s are semidegrees, it follows that $\scrC$ is indeed a cone. At first we show that $\scrC \subseteq \nef(\bar X)$. Pick positive integers $m_1, \ldots, m_k$ such that $D := \sum_{j=1}^k m_j [\tilde C_j] \in \scrC$. Fix $i$, $1 \leq i \leq k$, and $f_i \in \kk[X] \setminus \kk$ such that $\delta_i(f_i)/m_i \geq \delta_j(f_i)/m_j$ for all $j$, $1 \leq j \leq k$. The coefficient of $[\tilde C_j]$ in $D_i := \delta_i(f_i)D + \Div_{\bar X}(f_i^{m_i})$ is zero for $j = i$ and non-negative for $j \neq i$, $1 \leq j \leq k$. It follows that $(D, \tilde C_i) = (D_i, \tilde C_i)/\delta_i(f_i) \geq 0$ and therefore $D$ is nef, as required. \\

For the opposite inclusion, pick positive rational numbers $m_1, \ldots, m_k$ such that $D := \sum_{j=1}^k m_j [\tilde C_j]$ is in the {\em interior} of $\nef(\bar X)$. Then $D$ is ample, and therefore, replacing $D$ by some of its multiple if necessary we may assume that each $m_i$ is an integer and $D$ is very ample. Since $\sheaf_{\bar X}(D)$ is globally generated, for each $i$, $1 \leq i \leq k$, there exists $f_i \in \kk[X]$ such that $D_i := [\Div_{\bar X} (f)] + D \geq 0$ and the coefficient of $[\tilde C_i]$ in $D_i$ is zero. Then $\delta_i(f_i) = m_i$ and $\delta_j(f_i) \leq m_j$ for all $j \neq i$, so that the interior of $\nef(\bar X)$ is contained in $\scrC$. Since $\scrC$ is closed, we see that $\nef(\bar X) \subseteq \scrC$.
\end{proof}

Using Lemma \ref{nef-lemma} we give a positive answer to Question \ref{finite-question} in a special case. Recall that a degree-like function $\delta$ on a $\kk$-algebra is called {\em projective} if $\delta$ is finitely generated and $\delta(f) > 0$ for all $f \in A\setminus \kk$. 

\begin{cor}\label{finite-corollary}
Let $X$ be an affine surface and $\tilde \delta_1, \ldots, \tilde \delta_k$ be projective semidegrees on $\kk[X]$. Then $\tilde \delta := \max\{\tilde \delta_1, \ldots, \tilde \delta_k\}$ is finitely generated (and therefore also projective).  
\end{cor}

\begin{rem}
Not all affine varieties $X$ admit projective semidegrees on $\kk[X]$; a necessary condition for this to hold is the non-existence of any invertible regular functions on $X$ other than (non-zero) constants. 
\end{rem}

\begin{proof}
W.l.o.g.\ we may assume that $\tilde \delta = \max\{\tilde \delta_1, \ldots, \tilde \delta_k\}$ is the minimal presentation of $\tilde \delta$. At first assume that $X$ is normal and let $\bar X$ be the normalization of the closure in $\bar X_1 \times \cdots \times \bar X_k$ of the image of $X$ under the diagonal mapping. For each $j$, $1 \leq j \leq k$, define $d_j := \gcd\{\tilde \delta(f): f \in \kk[X]\}$, $\delta_j := \tilde \delta_j/d_j$ and $m_j := 1/d_j$. Lemma \ref{nef-lemma} then implies that $D := \sum_{j=1}^k m_j \tilde C_j$ is in the {\em interior} of the nef cone of $\bar X$, and therefore $D$ is ample. Since $D = \divinfinity{\tilde \delta}$, it follows that $\tilde \delta$ is finitely generated, as required to prove the corollary in the case that $X$ is normal. The general case follows (due to Remark \ref{normalizing-remark}) via applying the corollary to the normalization of $X$. 
\end{proof}

\bibliographystyle{alpha}
\bibliography{../../utilities/bibi}


\end{document}